\documentclass{article}
\usepackage[utf8]{inputenc}

\usepackage{fullpage}

\usepackage{amsmath}

\usepackage{amssymb}
\usepackage[mathscr]{eucal}
\usepackage{bbm}
\usepackage{etoolbox}

\newcommand{\RR}{\mathbb{R}}

\newcommand{\Aug}[1]{\widetilde{#1}}

\let\Vec\relax
\newcommand{\Vec}[1]{\boldsymbol{\mathrm{#1}}}

\newcommand{\Mat}[1]{\boldsymbol{\mathrm{#1}}}
\newcommand{\AugMat}[1]{\Aug{\Mat{#1}}{}}

\newcommand{\Tns}[1]{\boldsymbol{\mathscr{#1}}}
\newcommand{\AugTns}[1]{\Aug{\Tns{#1}}}

\newcommand{\TnsMat}[2]{\Mat{#1}_{(#2)}}
\newcommand{\AugTnsMat}[2]{\AugMat{#1}_{(#2)}}

\newcommand{\PartCore}[2]{\Tns{#1}^\ast_{#2}}

\newcommand{\PartCoreMat}[3]{\Mat{#1}^\ast_{#2,(#3)}}
\newcommand{\PartCoreMatT}[3]{\Mat{#1}^{\ast,\top}_{#2,(#3)}}

\newcommand{\GroupTnsMat}[3]{\Mat{#1}_{#2,(#3)}}

\newcommand{\Norm}[2][]{\lVert#2\rVert\ifstrempty{#1}{}{_{#1}}}

\usepackage{graphicx}
\usepackage{subcaption}

\usepackage{enumitem}

\usepackage{booktabs}
\usepackage{multirow}

\usepackage{algpseudocode}

\usepackage{amsthm}
\usepackage{algorithm}
\usepackage[capitalise]{cleveref}
\usepackage{url}

\newtheorem{proposition}{Proposition}

\theoremstyle{definition}
\newtheorem{assumption}[proposition]{Assumption}

\begin{document}

\title{Efficient Computation of Tucker Decomposition for Streaming Scientific Data Compression%
  \thanks{%
    This work was funded through the Advanced Scientific Computing Research (ASCR) program by the Office of Science (SC), U.S. Department of Energy.
    Sandia National Laboratories is a multimission laboratory managed and operated by National Technology \& Engineering Solutions of Sandia, LLC, a wholly owned subsidiary of Honeywell International Inc., for the U.S. Department of Energy's National Nuclear Security Administration under contract DE-NA0003525.
    This paper describes objective technical results and analysis. Any subjective views or opinions that might be expressed in the paper do not necessarily represent the views of the U.S. Department of Energy or the United States Government.
    This article has been authored by an employee of National Technology \& Engineering Solutions of Sandia, LLC under Contract No. DE-NA0003525 with the U.S. Department of Energy (DOE).
    The employee owns all right, title and interest in and to the article and is solely responsible for its contents.
    The United States Government retains and the publisher, by accepting the article for publication, acknowledges that the United States Government retains a non-exclusive, paid-up, irrevocable, world-wide license to publish or reproduce the published form of this article or allow others to do so, for United States Government purposes.
    The DOE will provide public access to these results of federally sponsored research in accordance with the DOE Public Access Plan \texttt{https://www.energy.gov/downloads/doe-public-access-plan}.
  }
}
\author{
  Saibal De\thanks{Sandia National Laboratories, Livermore, California, USA 94550 (\texttt{sde@sandia.gov}, \texttt{hnkolla@sandia.gov}).}
  \and
  Zitong Li\thanks{University of California, Irvine, California, USA 92697 (\texttt{zitongl5@uci.edu}).}
  \and
  Hemanth Kolla\footnotemark[2]
  \and
  Eric T. Phipps\thanks{Sandia National Laboratories, Albuquerque, New Mexico, USA 87123 (\texttt{etphipp@sandia.gov}).}
}
\date{}

\maketitle

\begin{abstract}

The Tucker decomposition, an extension of singular value decomposition for higher-order tensors, is a useful tool in analysis and compression of large-scale scientific data.
While it has been studied extensively for static datasets, there are relatively few works addressing the computation of the Tucker factorization of streaming data tensors.
In this paper we propose a new streaming Tucker algorithm tailored for scientific data, specifically for the case of a data tensor whose size increases along a single streaming mode that can grow indefinitely, which is typical of time-stepping scientific applications.
At any point of this growth, we seek to compute the Tucker decomposition of the data generated thus far, without requiring storing the past tensor slices in memory.
Our algorithm accomplishes this by starting with an initial Tucker decomposition and updating its components--the core tensor and factor matrices--with each new tensor slice as it becomes available, while satisfying a user-specified threshold of norm error.
We present an implementation within the TuckerMPI software framework, and apply it to synthetic and combustion simulation datasets.
By comparing against the standard (batch) decomposition algorithm we show that our streaming algorithm provides significant improvements in memory usage.
If the tensor rank stops growing along the streaming mode, the streaming algorithm also incurs less computational time compared to the batch algorithm.

\end{abstract}


\section{Introduction}
\label{sec:intro}

The push to exascale computing has led to tremendous advancements in the ability of scientific computing to generate high-fidelity data representing complex physics across wide ranges of spatial and temporal scales.
Similarly, modern sensing and data acquisition technologies have made real-time measurement of these same physical processes practical.
However, advances in data storage technologies have been unable to keep pace, making traditional offline analysis workflows impractical for scientific discovery.
Reduced representations that retain essential information, while guaranteeing a user-specified accuracy measure, have become an imperative for most scientific applications.

Conventional techniques for data compression (e.g. image compression) are not effective at achieving the levels of reduction required for scientific data, and various techniques specifically tailored for scientific data compression have been developed recently.
One of these, the Tucker decomposition, is the focus of our study.
The Tucker decomposition is one form of tensor decompositions which are akin to matrix factorizations extended to higher dimensions.
Tensors, or multidimensional data arrays, are a natural means of representing multivariate, high-dimensional data and tensor decompositions offer mathematically compact, computationally inexpensive, and inherently explainable models of scientific data.
Tucker decomposition has been shown to be quite efficient at scientific data compression, achieving orders-of-magnitude reduction for modest loss of accuracy \cite{BallardKK2020}.

The algorithm and implementation of a distributed-memory Tucker decomposition presented by Ballard {\it et al.} \cite{BallardKK2020} has very good parallel performance and scalability.
Nonetheless, it is designed to perform the decomposition in an off-line mode, {\it i.e.} perform the decomposition and compression of the data after it has already been saved to persistent file storage. 
Ballard {\it et al.} report that the memory footprint of their algorithm is about 3$\times$ the size of the original data, i.e. running the Tucker compression algorithm requires at least three times as much memory as the input data (see Sec. 9.4.1 of \cite{BallardKK2020}).
As datasets produced by scientific simulations get ever larger, this constraint can make it prohibitive to run Tucker compression in an off-line mode.
In practice, the need for compression is to avoid having to write large volumes of data to file storage in the first place.
To accomplish this the Tucker decomposition has to be performed on the data as it is being generated, which is challenging since most scientific applications have a streaming aspect to their data---a common scenario being time stepping on a grid representing a spatial domain---and the tensor is constantly growing along one or more dimensions.
This motivates our study and our objective is \emph{to develop a computationally efficient algorithm for performing Tucker decomposition on streaming scientific data.}
We address this by:
\begin{itemize}
  \item
    Proposing an algorithm for incremental update of the Tucker model by admitting new samples of streaming scientific data as tensor slices, while \emph{provably} satisfying error constraints of the conventional Tucker algorithm.
  \item
    Drawing motivation from ideas underpinning incremental singular value decomposition algorithms and delineating similarities to these ideas for the streaming and non-streaming tensor modes in our algorithm.
  \item
    Discussing the scalability of computational kernels and presenting performance improvements, relative to the batch (non-streaming) Tucker algorithm, for different data sizes.
\end{itemize}
In a prior work \cite{DeKoMePhRi2024}, we have analyzed the parallel scalability and portability of an earlier iteration of this streaming Tucker algorithm.
In this work, we primarily focus on the theoretical properties, specifically the approximation error bounds, of the algorithm in a serial setting.

The outline of our paper is as follows.
\Cref{sec:background} presents a background on tensor decompositions and the notations, definitions, and algorithms pertinent to Tucker factorization.
\Cref{sec:related_work} summarizes prior and related work on tensor decompositions for streaming data.
\Cref{sec:algorithm} contains the primary contribution of this paper, an algorithm for Tucker decomposition of streaming data.
We present the main building blocks of our algorithm, including a brief review of incremental singular value decomposition, alongside their theoretical performances.
\Cref{sec:results} presents some experimental results; using a synthetic dataset and a simulation dataset of turbulent combustion, we assess the accuracy and performance improvements (flops, memory usage) of the streaming algorithm over the batch version.


\section{Background}
\label{sec:background}

In this section, we first establish the notations we use in the remainder of this paper.
We then briefly review the definition of the Tucker decomposition and how it is computed from tensor data.
We conclude the section by providing an overview of streaming tensor data setup.

\subsection{Notations and Definitions}
\label{sec:notations}

We use bold lowercase letters (e.g., $\Vec{b}$) to denote vectors, and bold uppercase letters (e.g., $\Mat{U}$) to denote matrices.
The symbols $\Mat{0}$ and $\Mat{I}$ are reserved to denote zero vectors/matrices and the identity matrices; their sizes are sometimes specified through subscripts for clarity.
A $m \times n$ matrix $\Mat{U}$ is called orthogonal if $m \geq n$ and $\Mat{U}^\top \Mat{U} = \Mat{I}$, and row-orthogonal if its transpose is orthogonal (equivalently, $m \leq n$ and $\Mat{U} \Mat{U}^\top = \Mat{I}$).

Tensors are simply multi-way (multi-dimensional) arrays, with each dimension referred to as a \emph{mode}.
In this context, a matrix is a 2-way tensor with mode-1 being the row dimension, and mode-2 being the column dimension.
We use uppercase letters in script font (e.g.\ $\Tns{X}$) to denote tensors.
If $\Tns{X}$ is a $d$-way tensor of size $n_1 \times \dots \times n_d$, then we write $\Tns{X} \in \RR^{n_1 \times \cdots \times n_d}$.
To specify an element or a subset of the elements of a tensor, we borrow the subscripting and slicing notations from MATLAB.
For example, the first element of a 3-way (or 3D) tensor $\Tns{X}$ is $\Tns{X}(1,1,1)$, the first mode-1 fiber of $\Tns{X}$ is $\Tns{X}(:,1,1)$, and the first three frontal slices of $\Tns{X}$, taken together, is $\Tns{X}(:,:,1:3)$.

Two tensor operations are key to the Tucker algorithm.
The first is the matricization of a tensor along any mode that converts the tensor into a matrix.
The mode-$k$ matricization a tensor $\Tns{X}$ of size $n_1 \times \cdots \times n_d$ is the matrix of size $n_k \times n_1 \cdots n_{k - 1} n_{k + 1} \cdots n_d$, denoted as $\TnsMat{X}{k}$.
The $k$-th mode of the tensor corresponds to the row dimension of the matrization, and the remaining modes are combined to form the column dimension; entries of this matrix using the index notation are given by
\begin{align*}
  \TnsMat{X}{k}(i_k, \overline{i_1, \ldots, i_{k - 1}, i_{k + 1}, \ldots, i_d}) &= \Tns{X}(i_1, \ldots, i_d), \\
  \overline{i_1, \ldots, i_{k - 1}, i_{k + 1}, \ldots, i_d} &= 1 + \sum_{\substack{\ell = 1\\\ell \neq k}}^d (i_\ell - 1) \prod_{\substack{\ell' = 1\\\ell' \neq k}}^{\ell - 1} n_{\ell'}.
\end{align*}

The second operation is the tensor-times-matrix (TTM) product along a specified mode of the tensor.
The TTM of a tensor $\Tns{X} \in \RR^{n_1 \times \cdots \times n_d}$ and a matrix $\Mat{A} \in \RR^{m \times n_k}$ along mode-$k$ produces a tensor $\Tns{Y} \in \RR^{n_1 \times \cdots \times n_{k - 1} \times m \times n_{k + 1} \times \cdots \times n_d}$.
This operation, denoted by $\Tns{Y} = \Tns{X} \times_k \Mat{A}$, modifies the input tensor by multiplying all of its mode-$k$ fibers by the matrix; the entries of $\Tns{X}$ and $\Tns{Y}$ are related via the matrix product
\begin{equation*}
  \TnsMat{Y}{k} = \Mat{A} \TnsMat{X}{k}.
\end{equation*}
TTM operation along different modes commute: $\Tns{X} \times_k \Mat{A} \times_{k'} \Mat{B} = \Tns{X} \times_{k'} \Mat{B} \times_k \Mat{A}$ for any $k \neq k'$ and appropriately sized matrices $\Mat{A}$ and $\Mat{B}$.
If two $d$-dimensional tensors $\Tns{X} \in \RR^{m_1 \times \cdots \times m_d}$ and $\Tns{Y} \in \RR^{n_1 \times \cdots \times n_d}$ are related via a series of TTMs,
\begin{equation*}
  \Tns{Y} = \Tns{X} \times_1 \Mat{A}_1 \cdots \times_d \Mat{A}_d,
\end{equation*}
then their matricizations are also related as follows,
\begin{equation}
  \label{eq:tucker_matrix}
  \TnsMat{Y}{k} = \Mat{A}_k \TnsMat{X}{k} (\Mat{A}_d \otimes \cdots \otimes \Mat{A}_{k + 1} \otimes \Mat{A}_{k - 1} \otimes \cdots \otimes \Mat{A}_1)^\top, \quad 1 \leq k \leq d,
\end{equation}
where ``$\otimes$'' represents the matrix Kronecker product.
If the TTMs are not applied across all tensor modes, e.g.,
\begin{equation*}
  \Tns{Y} = \Tns{X} \times_1 \Mat{A}_1 \cdots \times_k \Mat{A}_k, \quad 1 \leq k < d,
\end{equation*}
then we necessarily have $m_\ell = n_\ell$ for $k + 1 \leq \ell \leq d$.
In this case, we can relate the matricizations of $\Tns{X}$ and $\Tns{Y}$ by filling in the missing TTM slots with identity matrices:
\begin{equation*}
  \Tns{Y} = \Tns{X} \times_1 \Mat{A}_1 \cdots \times_k \Mat{A}_k \times_{k + 1} \Mat{A}_{k + 1} \cdots \times_d \Mat{A}_d, \quad \Mat{A}_\ell = \Mat{I}_{n_\ell} \text{ for } k + 1 \leq \ell \leq d.
\end{equation*}

\subsection{Tucker Decomposition and Data Compression}
\label{sec:tucker}

The Tucker decomposition \cite{Tucker1966} factorizes a $n_1 \times \cdots \times n_d$ tensor $\Tns{X}$ in the form
\begin{equation}
  \label{eq:tucker_def}
  \Tns{X} = \Tns{C} \times_1 \Mat{U}_1 \cdots \times_d \Mat{U}_d,
\end{equation}
where $\Tns{C}$ is called the `core' tensor and $\Mat{U}_k$ are orthogonal `factor' matrices.
For a full rank decomposition, the core tensor is of size $n_1 \times \cdots \times n_d$ and each factor matrix $\Mat{U}_k$ has size $n_k \times n_k$ for $1 \leq k \leq d$.
Using the TTM matricization property from \cref{eq:tucker_matrix}, we can rewrite this Tucker decomposition in an equivalent form,
\begin{equation}
  \TnsMat{X}{k} = \Mat{U}_k \TnsMat{C}{k} (\Mat{U}_d \otimes \cdots \otimes \Mat{U}_{k + 1} \otimes \Mat{U}_{k - 1} \otimes \cdots \otimes \Mat{U}_1)^\top, \quad 1 \leq k \leq d.
\end{equation}
Aside from the core matricization $\TnsMat{C}{k}$ not necessarily being diagonal, this form is very similar to the matrix singular value decomposition (SVD).
In fact, each factor matrix $\Mat{U}_k$ is the orthonormal basis for the mode-$k$ matricization of the tensor, comprising of the left singular vectors of $\TnsMat{X}{k}$.
This way, the Tucker decomposition is a higher-order generalization of singular value decomposition (HOSVD) \cite{LathauwerMV2000}.
The core tensor is constructed from the data tensor by
\begin{equation*}
  \Tns{C} = \Tns{X}\times_1 \Mat{U}_1^\top \cdots \times_d \Mat{U}_d^\top.
\end{equation*}

\begin{figure}[t]
  \centering
  \includegraphics[width=0.5\linewidth]{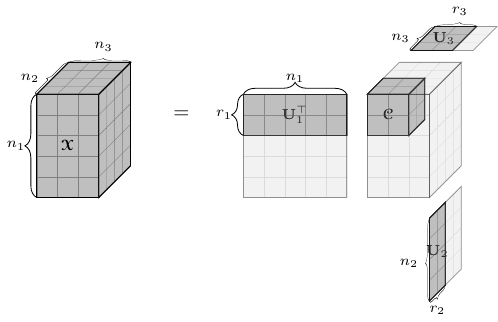}
  \caption{Tucker decomposition of a $3$-way tensor $\Tns{X}$ with size $n_1\times n_2 \times n_3$ and Tucker ranks $(r_1, r_2, r_3)$, consisting of a core tensor $\Tns{C}$ of size $r_1 \times r_2 \times r_3$ and orthgonal factor matrices $\Mat{U}_k \in \RR^{n_k \times r_k}$ for $1 \leq k \leq 3$.}
  \label{fig:tucker_decomp}
\end{figure}

If $\Tns{X}$ has low rank structure then the HOSVD results in a core tensor $\Tns{C}$ of size $r_1 \times \cdots \times r_d$, and factor matrix $\Mat{U}_k$ of size $n_k \times r_k$, with $r_k \leq n_k$; this is depicted in \cref{fig:tucker_decomp}.
The sizes $r_1, \ldots, r_d$ of the core tensor are called the Tucker ranks of the tensor.
Clearly, when the ranks are far smaller than the mode sizes $(r_k \ll n_k)$, we can achieve significant savings in storing the tensor in Tucker format;
the corresponding compression ratio is given by
\begin{equation*}
  \rho = \frac{n_1 \cdots n_d}{r_1 \cdots r_d + n_1 r_1 + \cdots + n_d r_d}.
\end{equation*}

\subsection{Sequentially Truncated Higher-Order SVD}

More often instead of an exact low-rank representation as in \cref{eq:tucker_def}, we seek an approximation to the data tensor.
We can compute such a low-rank representation via a series of truncated SVDs; this algorithm is known as the truncated higher order SVD (T-HOSVD) \cite{VVM12}.
It loops over the modes $k = 1, \ldots, d$ of a data tensor $\Tns{X}$, and computes the principal left singular vectors $\Mat{U}_k$ of the corresponding unfolding matrices $\TnsMat{X}{k}$.
Then the core $\Tns{C}$ is obtained by projecting the data tensor along these singular vector spans: $\Tns{C} = \Tns{X} \times_1 \Mat{U}_1^\top \cdots \times \Mat{U}_d^\top$.
The resulting Tucker approximation satisfies
\begin{equation}
  \label{eq:thosvd_err}
  \Norm[F]{\Tns{X} - \Tns{C} \times_1 \Mat{U}_1 \cdots \times_d \Mat{U}_d} \leq \tau \Norm[F]{\Tns{X}},
\end{equation}
where $\Norm[F]{\cdot}$ denotes the Frobenius norm and $\tau > 0$ is the error tolerance specified to the algorithm.

Vannieuwenhoven {\it et al.} \cite{VVM12} present an important refinement of T-HOSVD, where the computation is much more efficient and the error guarantee, \cref{eq:thosvd_err}, is still preserved.
This algorithm works with ``partial cores'', constructed by sequentially projecting the data tensor along already computed principal mode subspaces:
\begin{equation*}
  \PartCore{C}{k}
  = \Tns{X} \times_1 \Mat{U}_1^\top \cdots \times_k \Mat{U}_k^\top
  \in \RR^{r_1 \times \cdots \times r_k \times n_{k + 1} \times \cdots \times n_d},
\end{equation*}
with $\PartCore{C}{0} = \Tns{X}$ and $\PartCore{C}{d} = \Tns{C}$.
At each step of the algorithm, the partial core constructed in the last step is used to construct the principal subspace of the current tensor mode.
Note that in this sequentially truncated higher-order SVD (ST-HOSVD) variant, the order in which the modes are processed is important; this is in contrast with the T-HOSVD algorithm where modes can be processed in any order and produce the same final decomposition.
In this work, to avoid introducing additional notations, we assume the modes are processed in the order $k = 1, \ldots, d$.
During the computation, the partial cores and the tensor reconstructions generally satisfy the error bounds
\begin{equation}
  \label{eq:partial_core_error}
  \Norm[F]{\PartCore{C}{k - 1} - \PartCore{C}{k} \times \Mat{U}_k}
  \leq \frac{\tau}{\sqrt{d}} \Norm[F]{\Tns{X}},
  \quad 1 \leq k \leq d.
\end{equation}
The final error bound in \cref{eq:thosvd_err} is obtained by combining these error bounds over all tensor modes:
\begin{equation*}
  \Norm[F]{\Tns{X} - \Tns{C} \times_1 \Mat{U}_1 \cdots \times_d \Mat{U}_d}^2 = \sum_{k = 1}^d \Norm[F]{\PartCore{C}{k - 1} - \PartCore{C}{k} \times_k \Mat{U}_k}^2 \leq \sum_{k = 1}^d \frac{\tau^2}{d} \Norm[F]{\Tns{X}}^2 = \tau^2 \Norm[F]{\Tns{X}}^2.
\end{equation*}

\begin{algorithm}[t]
  \caption{Truncated Gram SVD}
  \label{alg:tsvd}
  \begin{algorithmic}[1]
    \Function{TruncatedSVD}{$\Mat{A}$, $\epsilon$} \Comment{Matrix $\Mat{A} \in \RR^{m \times n}$, absolute tolerance $\epsilon \geq 0$}
      \State {Compute Gram matrix: $\Mat{G} \gets \Mat{A} \Mat{A}^\top$}
      \State {Eigendecompose: $\Vec{\lambda}, \Mat{P} \gets \Call{eig}{\Mat{G}}$}
      \State {Determine rank: $r \gets \min \{\ell : \lambda_{\ell + 1} + \cdots + \lambda_m \leq \epsilon^2\}$}
      \State {Truncated left singular vectors: $\Mat{U} \gets \Mat{P}(:, 1 : r)$}
      \State {Truncated singular values: $\Mat{S} \gets \Call{diag}{\Vec{\lambda}(1 : r)}$} \Comment{Optional}
      \State {Truncated right singular vectors: $\Mat{V} \gets \Mat{A}^\top \Mat{U} \Mat{S}^{-1}$} \Comment{Optional}
      \State \Return{$r$, $\Mat{U}$, $\Mat{S}$, $\Mat{V}$}
    \EndFunction \Comment{$\Norm[F]{\Mat{A} - \Mat{U} \Mat{S} \Mat{V}^\top} \leq \epsilon$}
  \end{algorithmic}
\end{algorithm}

\begin{algorithm}[t]
  \caption{Sequentially truncated higher order SVD (ST-HOSVD)}
  \label{alg:sthosvd}
  \begin{algorithmic}[1]
    \Function{ST-HOSVD}{$\Tns{X}$, $\tau$} \Comment{Tensor $\Tns{X} \in \RR^{n_1 \times \cdots \times n_d}$, relative tolerance $\tau > 0$}
      \State{Compute truncation: $\epsilon \gets \tau \Norm[F]{\Tns{X}} / \sqrt{d}$}
      \State{Initialize partial core: $\PartCore{C}{0} \gets \Tns{X} $}
      \For{$k = 1$ to $d$}
        \State{Determine truncated mode basis: $r_k, \Mat{U}_k \gets \Call{TruncatedSVD}{\PartCoreMat{C}{k - 1}{k}, \epsilon}$}
        \State{Update partial core: $\PartCore{C}{k} \gets \PartCore{C}{k - 1} \times_k \Mat{U}_k^\top$}
      \EndFor
      \State{Compute core: $\Tns{C} \gets \PartCore{C}{d}$}
      \State \Return{$\Tns{C}$, $\Mat{U}_1, \ldots, \Mat{U}_d$}
    \EndFunction \Comment{$\Norm[F]{\Tns{X} - \Tns{C} \times_1 \Mat{U}_1 \cdots \times_d \Mat{U}_d} \leq \tau \Norm[F]{\Tns{X}}$}
  \end{algorithmic}
\end{algorithm}

In \cite{AustinBK2016,BallardKK2020}, the authors propose a further modification; they compute the SVD of the partial core unfolding $\PartCoreMat{C}{k - 1}{k}$ of size $n_k \times r_1 \cdots r_{k - 1} n_{k + 1} \cdots n_d$ through eigendecomposition of the much smaller Gram matrix of size $n_k \times n_k$.
The loss of precision of this Gram-SVD approach is generally acceptable since the tolerance $\tau$ is typically higher than the square root of machine precision.
The steps of this ST-HOSVD method are outlined in \cref{alg:sthosvd}.

\subsection{Problem Definition: Streaming Scientific Data}

In recent years, numerous methods for computing tensor decompositions of streaming data have been developed.
Unfortunately, there is no single well-defined formulation of the streaming problem, so the assumptions made by these methods vary.
In this work, we make the following assumptions:
Data arrives in discrete batches, which we call time steps indexed by $t = 1, 2, \ldots$, and batches will be processed sequentially.
Furthermore, each batch is assumed to be a $(d - 1)$-way tensor, which we will call $\Tns{Y} \in \RR^{n_1 \times \cdots \times n_{d - 1}}$, and the dimensions $n_1, \ldots, n_{d - 1}$ of each batch fixed throughout time.
Finally, we assume the number of time steps is unbounded making it impossible to store all observed data, or revisit sufficiently old data.
These assumptions are clearly satisfied in the context of scientific computing, where each batch corresponds to a time step in a scientific simulation, but most other data analysis problems can be reformulated to fit these assumptions.
Within this context, there are multiple equivalent ways of viewing the streaming problem.
The most common, and the approach considered here, is to stack the batches $\Tns{Y}$ along a new ``temporal'' mode $d$, resulting in a $d$-way tensor $\Tns{X} \in \RR^{n_1 \times \cdots \times n_{d - 1} \times n_d}$ where $n_d$ is the number of batches.


\section{Related Work}
\label{sec:related_work}

While much prior work on streaming canonical polyadic (CP) tensor decompositions exists, e.g., \cite{NiSi09, MaMaGi15, SmHuSiKa18, ZhViBaJi16, Ka16, SoHuGeCa17, GuPaPa18, PaGuPa19, LeBaHeEz18, PhCi11a, VaVeLa17}, streaming methods for Tucker decompositions are much less well developed.
Early work, called Incremental Tensor Analysis~\cite{Sun2008}, proposed several algorithms for computing Tucker decompositions of a stream of tensors.
Several later works~\cite{Sobral2014,Ma2009,Hu2011} improved the performance of these approaches by leveraging incremental SVD techniques and modifying them to better suit computer vision use cases.
Regardless, we find these approaches less appealing since a new core tensor is computed for each streamed tensor, and there is no compression over time.
Conversely,~\cite{Shi2015} only updates the temporal mode factor matrix of a tensor stream in support of anomaly detection, and never computes a full Tucker decomposition.
Recently two approaches \cite{MaBe18,Sun2020} relying on sketching have been developed for Tucker decompositions that only require a single pass over the tensor data, and therefore can be used to decompose streaming data.
However, in both cases the methods cannot form the Tucker decomposition until all data has been processed, and therefore are inappropriate for analysis of unbounded data streams.
Additionally~\cite{Xiao2018} considers potential growth in all modes over time, not just a temporal mode, but replaces the usual SVD calculations with matrix inverses, bringing into question the accuracy and scalability of the method.
Moreover, the algorithm in \cite{Xiao2018} only addresses the case of an exact Tucker decomposition, not a truncated or sequentially truncated decomposition where the model is an approximation of the original tensor, which is necessary for lossy compression.
Finally, \cite{JaKa2023} introduces the D-Tucker and D-TuckerO algorithms for constructing static and online Tucker factorization using randomized SVD of the 2D slices of the data tensor; however this approach is fundamentally different from the HOSVD and ST-HOSVD algorithms we consider in this paper.

Given the close relationship between Tucker decompositions and matrix SVD, several of the approaches above formulate a streaming Tucker algorithm in terms of incremental SVD, for which several methods have been developed, e.g.,~\cite{Brand02,Brand06, Zhang22,Ross08,Cheng18}.
Our approach is similar and applies an incremental SVD approach inspired by Brand~\cite{Brand02} to the streaming (temporal) mode, but pursues a different approach for updating the non-streaming modes.


\section{Algorithm for Streaming ST-HOSVD}
\label{sec:algorithm}

We now present our algorithm for applying streaming updates to an existing Tucker factorization.
In the following, we first briefly review the key concepts of incremental SVD update for streaming matrix factorization, then develop the streaming Tucker update algorithm.

\subsection{Incremental SVD}
\label{sec:incremental-svd}

Let $\Mat{A} \approx \Mat{U} \Mat{S} \Mat{V}^\top$ be the rank-$r$ approximate SVD of data matrix $\Mat{A} \in \RR^{m \times n}$ with orthogonal $\Mat{U} \in \RR^{m \times r}$, $\Mat{V} \in \RR^{n \times r}$ and diagonal $\Mat{S} \in \RR^{r \times r}$.
As a new row $\Vec{b} \in \RR^n$ is added to this matrix, we want to update the SVD factors $\Mat{U}$, $\Mat{S}$ and $\Mat{V}$ to construct an approximation of the augmented matrix
\begin{equation*}
  \AugMat{A} =
  \begin{bmatrix}
    \Mat{A} \\
    \Vec{b}^\top
  \end{bmatrix}
  \in \RR^{(m + 1) \times n}.
\end{equation*}

We follow the basic version of the incremental SVD algorithms presented in \cite{Brand02,Cheng18}.
First split the new row into orthogonal components w.r.t.\ the column span of the right singular vectors $\Mat{V}$:
\begin{equation*}
  \Vec{p} = \Mat{V}^\top \Vec{b} \in \RR^{r},
  \quad
  \Vec{e} = \Vec{b} - \Mat{V} \Vec{p} \in \RR^n
  \implies
  \Vec{b} = \Mat{V} \Vec{p} + \Vec{e}.
\end{equation*}
If the orthogonal complement $\Vec{e}$ is approximately zero, then observe that
\begin{equation}
  \label{eq:isvd_main_1}
  \AugMat{A}
  =
  \begin{bmatrix}
    \Mat{A} \\
    \Vec{b}^\top
  \end{bmatrix}
  \approx
  \begin{bmatrix}
    \Mat{U} \Mat{S} \Mat{V}^\top \\
    \Vec{p}^\top \Mat{V}^\top
  \end{bmatrix}
  =
  \begin{bmatrix}
    \Mat{U} &   \\
            & 1
  \end{bmatrix}
  \begin{bmatrix}
    \Mat{S} \\
    \Vec{p}^\top
  \end{bmatrix}
  \Mat{V}^\top.
\end{equation}
Otherwise, if $\Vec{e}$ is non-zero, then construct a new basis vector for the row-span of the augmented data matrix,
\begin{equation*}
  q = \Norm{\Vec{e}},
  \quad
  \Vec{v} = q^{-1} \Vec{e}
  \implies
  \Vec{b} = \Mat{V} \Vec{p} + q \Vec{v},
\end{equation*}
and observe that
\begin{equation}
  \label{eq:isvd_main_2}
  \AugMat{A} =
  \begin{bmatrix}
    \Mat{A} \\
    \Vec{b}^\top
  \end{bmatrix} \approx
  \begin{bmatrix}
    \Mat{U} \Mat{S} \Mat{V}^\top \\
    \Vec{p}^\top \Mat{V}^\top + q \Vec{v}^\top
  \end{bmatrix}
  =
  \begin{bmatrix}
    \Mat{U} &   \\
            & 1
  \end{bmatrix}
  \begin{bmatrix}
    \Mat{S}      &   \\
    \Vec{p}^\top & q
  \end{bmatrix}
  \begin{bmatrix}
    \Mat{V} & \Vec{v}
  \end{bmatrix}^\top.
\end{equation}
In either case, we have factorized the augmented matrix into a product of the form $\AugMat{A} \approx \Mat{U}' \Mat{S}' \Mat{V}^{\prime,\top}$.
It is easy to verify that $\Mat{U}'$ and $\Mat{V}'$ are orthogonal matrices, and we can see that $\Mat{S}'$ is nearly diagonal.
We compute the SVD of this middle matrix, $\Mat{S}' = \Mat{U}'' \AugMat{S} \Mat{V}^{\prime\prime,\top}$, to finally obtain SVD of the augmented matrix:
\begin{equation}
  \label{eq:isvd_final}
  \AugMat{A} \approx \AugMat{U} \AugMat{S} \AugMat{V}^\top,
  \quad
  \AugMat{U} = \Mat{U}' \Mat{U}'',
  \quad
  \AugMat{V} = \Mat{V}' \Mat{V}''.
\end{equation}

\begin{algorithm}[t]
  \caption{Incremental SVD}
  \label{alg:isvd}
  \begin{algorithmic}[1]
    \Function{ISVD}{$\Mat{U}$, $\Mat{S}$, $\Mat{V}$, $\Vec{b}$, $\epsilon$} \Comment{SVD $\Mat{U}$, $\Mat{S}$, $\Mat{V}$ of $\Mat{A}$, new row $\Vec{b}$, absolute tol.\ $\epsilon$}
      \State {$\Vec{p} \gets \Mat{V}^\top \Vec{b}$}
      \State {$\Vec{e} \gets \Vec{b} - \Mat{V} \Vec{p}$}
      \State {$q \gets \| \Vec{e} \|$}
      \If {$q \leq \epsilon$}
        \State {$\Mat{S}' \gets \begin{bmatrix} \Mat{S} \\ \Vec{p}^\top \end{bmatrix}$}
        \State {$\Mat{V}' \gets \Mat{V}$}
      \Else
        \State {$\Vec{v} \gets q^{-1} \Vec{e}$}
        \State {$\Mat{S}' \gets \begin{bmatrix} \Mat{S} & \\ \Vec{p}^\top & q \end{bmatrix}$}
        \State {$\Mat{V}' \gets \begin{bmatrix} \Mat{V} & \Vec{v} \end{bmatrix}$}
        \If {$\lvert \Mat{V}(:, 1)^\top \Vec{v} \rvert > \epsilon_\text{m}$} \Comment{Orthogonality check against some small $\epsilon_m$}
          \State {$\Mat{Q}, \Mat{R} \gets \Call{QR}{\Mat{V}'}$}
          \State {$\Mat{V}' \gets \Mat{Q}$}
          \State {$\Mat{S}' \gets \Mat{S}' \Mat{R}^\top$}
        \EndIf
      \EndIf
      \State {$\Mat{U}'', \AugMat{S}, \Mat{V}'' \gets \Call{TruncatedSVD}{\Mat{S}', 0}$} \Comment{Exact compact SVD}
      \State {$\AugMat{U} \gets \begin{bmatrix} \Mat{U} \Mat{U}''(1 : r, :) \\ \Mat{U}''(r + 1, :) \end{bmatrix}$} \Comment{$\AugMat{U} \gets \Mat{U}' \Mat{U}''$ without explicitly allocating $\Mat{U}'$}
      \State {$\AugMat{V} \gets \Mat{V}' \Mat{V}''$}
      \State \Return {$\AugMat{U}$, $\AugMat{S}$, $\AugMat{V}$}
    \EndFunction \Comment{$\Norm[F]{\AugMat{A} - \AugMat{U} \AugMat{S} \AugMat{V}{}^\top}^2 \leq \Norm[F]{\Mat{A} - \Mat{U} \Mat{S} \Mat{V}^\top}^2 + \epsilon^2$}
  \end{algorithmic}
\end{algorithm}

We quantify the approximations in the algorithm description in terms of the Frobenius norm.
Aside from the initial error in the approximation of $\Mat{A}$, potentially ignoring the orthogonal complement $\Vec{e}$ introduces additional errors to the approximation of the augmented matrix $\AugMat{A}$.
The precise characterization is as follows:

\begin{proposition}
  \label{prop:isvd-error-bound}
  The error in approximating the augmented matrix $\AugMat{A}$ using the updated factors from incremental SVD is given by
  \begin{equation*}
    \Norm[F]{\AugMat{A} - \AugMat{U} \AugMat{S} \AugMat{V}^\top}^2 =
    \begin{cases}
      \Norm[F]{\Mat{A} - \Mat{U} \Mat{S} \Mat{V}^\top} + \Norm{\Vec{e}}^2 & \text{if } \Vec{e} \text{ is ignored for being too small}, \\
      \Norm[F]{\Mat{A} - \Mat{U} \Mat{S} \Mat{V}^\top}                    & \text{if } \Vec{e} \text{ is used to create basis vector } \Vec{v}.
      \end{cases}
  \end{equation*}
  In particular, for a fixed $\epsilon > 0$, suppose we ignore $\Vec{e}$ if $\Norm{\Vec{e}} \leq \epsilon$.
  Then the incremental SVD error evolves as $\Norm[F]{\AugMat{A} - \AugMat{U} \AugMat{S} \AugMat{V}^\top}^2 \leq \Norm[F]{\Mat{A} - \Mat{U} \Mat{S} \Mat{V}^\top}^2 + \epsilon^2$.
\end{proposition}

\begin{proof}
  See \cref{app:isvd-error-bound}.
\end{proof}

A direct consequence of this error evolution is the following.
Suppose the initial approximation of $\Mat{A}$ has relative accuracy $\tau$, satisfying $\Norm[F]{\Mat{A} - \Mat{U} \Mat{S} \Mat{V}^\top} \leq \tau \Norm[F]{\Mat{A}}$.
Then by setting $\epsilon = \tau \Norm{\Vec{b}}$, we observe that
\begin{equation*}
  \Norm[F]{\AugMat{A} - \AugMat{U} \AugMat{S} \AugMat{V}^\top} \leq \sqrt{\tau^2 \Norm[F]{\Mat{A}}^2 + \tau^2 \Norm{\Vec{b}}^2} = \tau \sqrt{\Norm[F]{\Mat{A}}^2 + \Norm{\Vec{b}}^2} = \tau \Norm[F]{\AugMat{A}},
\end{equation*}
i.e., the relative accuracy of the updated factorization is maintained.

The addition of a new basis $\Vec{v}$ to the existing right singular vectors $\Mat{V}$ to form the updated $\Mat{V}'$ matrix may lead to loss of orthogonality due to floating point round off.
To resolve this, we occasionally use a QR factorization to ensure $\Mat{V}'$ remains orthogonal. \Cref{alg:isvd} details the key steps of this procedure.

The advantage of this incremental update is that, instead of computing the SVD of the augmented data matrix $\AugMat{A} \in \RR^{(m + 1) \times n}$ from scratch, we can use the factorization of the existing data matrix $\Mat{A}$ to assimilate the new data row $\Vec{b}$.
In the process, we compute the SVD of a \emph{smaller} matrix $\Mat{S}' \in \RR^{(r + 1) \times \Aug{r}}$, where $\Aug{r}$ is either $r$ or $r + 1$ depending on whether the orthogonal complement $\Vec{e}$ is sufficiently small.
This approach is suitable for the streaming data setup, where the data matrix is being generated one row at a time.
In addition, if the final data matrix is low rank, then incremental update is also computationally efficient.
Assuming that we are \emph{not} using a specialized algorithm to compute the SVD of $\Mat{S}'$ that exploits its almost diagonal structure, the cost of the incremental SVD is characterized as follows:

\begin{proposition}
  \label{prop:isvd-flop-count}
  The assimilation of a row to the rank-$r$ SVD of a $m \times n$ data matrix via incremental SVD uses $O(m r^2 + n r^2 + r^3)$ floating point operations and $O(m r + n r + r^2)$ extra storage.
  Assuming the rank $r$ does not change over the updates, the amortized cost of incremental SVD of a data matrix with final size $m \times n$ is $O(m^2 r^2 + m n r^2 + m r^3)$ operations and $O(m r + n r + r^2)$ storage.
\end{proposition}

\begin{proof}
  See \cref{app:isvd-flop-count}.
\end{proof}

Contrast this against the cost of computing the $r$-rank truncated Gram-SVD of a $m \times n$ matrix, which uses $O(m^3 + m^2 n + m n r + n r^2)$ operations and $O(m^2 + m r + n r + r^2)$ storage.
When the matrix has many more rows than columns---a natural assumption in the streaming setup---and the matrix rank remains small, then the leading order time and memory complexities of incremental SVD are $O(m^2)$ and $O(m)$.
These are an order of magnitude smaller than the leading order $O(m^3)$ time complexity and $O(m^2)$ memory complexity of batch SVD.

Readers should note that the algorithm described above is a basic version of the incremental SVD.
A drawback of this approach is the update of the singular vectors in \cref{eq:isvd_final}.
For example, in constructing the new left singular vectors, we multiply a $(m + 1) \times (r + 1)$ orthogonal matrix $\Mat{U}'$ with a smaller $(r + 1) \times \Aug{r}$ orthogonal matrix $\Mat{U}''$.
When the ranks are small (i.e.\ $m \gg r$), this may lead to loss-of-orthogonality after a large number of incremental updates due to numerical error accumulation.
In \cite{Brand02,Cheng18}, the authors present techniques to overcome this limitation; however, for the computational experiments presented here, the simplified algorithm works well enough.
We leave incorporating the improved algorithms to future work.

\subsection{Streaming Tensor Setup}

We now adapt the basic approach for incremental SVD, and develop a streaming Tucker update algorithm for higher-dimensional tensors.
Suppose $\Tns{X} \approx \Tns{C} \times_1 \Mat{U}_1 \cdots \times_d \Mat{U}_d$ is an approximate Tucker factorization of data tensor $\Tns{X}$.
Let us assume, for the moment, that this factorization was produced using ST-HOSVD (\cref{alg:sthosvd}) where the tensor modes were processed in order $k = 1, \ldots, d$, with the last mode being the streaming dimension%
\footnote{This assumption is for the sake of simplicity; the key requirement is that the streaming mode is processed last in the batch or streaming ST-HOSVD algorithms.}.
Then the core $\Tns{C}$ was constructed from the partial core $\PartCore{C}{d - 1}$ using truncated Gram SVD:
\begin{equation*}
  \Mat{G} = \PartCoreMat{C}{d - 1}{d} \PartCoreMatT{C}{d - 1}{d}, \quad \Mat{G} \approx \Mat{U}_d \Mat{S} \Mat{U}_d^\top, \quad \TnsMat{C}{d} \equiv \PartCoreMat{C}{d}{d} = \Mat{U}_d^\top \PartCoreMat{C}{d - 1}{d}.
\end{equation*}
Let $\PartCoreMat{C}{d - 1}{d} = \Mat{U}_d \Mat{S} \Mat{V}^\top + \Mat{U}_d' \Mat{S}' \Mat{V}^{\prime,\top}$ be the full SVD of the matricization $\PartCoreMat{C}{d - 1}{d}$, where the second term represents discarded singular values and singular vectors from the truncation.
Since the truncated SVD algorithm guarantees that $\Mat{U}_d$ is a orthogonal matrix, and that $\Mat{U}_d$ and $\Mat{U}_d'$ have mutually orthogonal column spans, it follows that
\begin{equation*}
  \TnsMat{C}{d}
  = \Mat{U}_d^\top \PartCoreMat{C}{d - 1}{d}
  = \Mat{U}_d^\top \Mat{U}_d \Mat{S} \Mat{V}^\top + \Mat{U}_d^\top \Mat{U}_d' \Mat{S}' \Mat{V}^{\prime,\top}
  = \Mat{I} \Mat{S} \Mat{V}^\top + \Mat{0} \Mat{S}' \Mat{V}^{\prime,\top}
  = \Mat{S} \Mat{V}^\top.
\end{equation*}
In other words, the core admits a further factorization into a diagonal matrix $\Mat{S}$ and another tensor $\Tns{V}$ of same size as $\Tns{C}$ with $\TnsMat{V}{d} = \Mat{V}^\top$ such that
\begin{equation}
  \label{eq:tucker-core-factorization}
  \Tns{C} = \Tns{V} \times_d \Mat{S} \iff \TnsMat{C}{d} = \Mat{S} \TnsMat{V}{d}.
\end{equation}
Rows of the matricization $\TnsMat{V}{d}$ correspond to the right singular vectors of the last truncated SVD in the ST-HOSVD procedure, and therefore the matricization is a row-orthogonal matrix.
This will become important as we develop our streaming ST-HOSVD update, so we add it to the list of assumptions we make about our current Tucker approximation:

\begin{assumption}
  \label{asm:streaming_setup}
  We are given an approximate Tucker factorization of data tensor $\Tns{X} \in \RR^{n_1 \times \cdots \times n_d}$.
  \begin{enumerate}[label=(\alph*)]
    \item
      The factorization consists of core $\Tns{C} \in \RR^{r_1 \times \cdots \times r_d}$ and factor matrices $\Mat{U}_k \in \RR^{n_k \times r_k}$ for $1 \leq k \leq d$.
    \item
      The approximation satisfies the error bound from \cref{eq:thosvd_err}, $\Norm[F]{\Tns{X} - \Tns{C} \times_1 \Mat{U}_1 \cdots \times_d \Mat{U}_d} \leq \tau \Norm[F]{\Tns{X}}$, where $\tau > 0$ is some fixed tolerance.
    \item
      The core $\Tns{C}$ admits a further factorization along the last tensor mode, $\Tns{C} = \Tns{V} \times_d \Mat{S}$, where $\TnsMat{V}{d}$ is row-orthogonal, and $\Mat{S}$ is diagonal with $\Mat{S}(i, i) > 0$ for all $1 \leq i \leq r_d$.
  \end{enumerate}
\end{assumption}

Given this initial Tucker factorization and a new tensor slice $\Tns{Y} \in \RR^{n_1 \times \cdots \times n_{d - 1} \times 1}$, we augment the data tensor along the last tensor mode,
\begin{equation}
  \label{eq:augmented_tensor}
  \AugTns{X}(i_1, \ldots, i_d) =
  \begin{cases}
    \Tns{X}(i_1, \ldots, i_d) & \text{if} \quad 1 \leq i_d \leq n_d, \\
    \Tns{Y}(i_1, \ldots, i_{d - 1}, 1) & \text{if} \quad i_d = n_d + 1,
  \end{cases}
\end{equation}
and seek its approximate Tucker factorization.
We design our streaming update algorithm in such way that:
\begin{enumerate}[label=(\alph*)]
\item
  The updated Tucker factorization consists of core $\AugTns{C} \in \RR^{\Aug{r}_1 \cdots \times \Aug{r}_d}$ and orthogonal factors $\AugMat{U}_k \in \RR^{n_k \times \Aug{r}_k}$ for $1 \leq k \leq d - 1$ and $\AugMat{U}_d \in \RR^{(n_d + 1) \times \Aug{r}_d}$.
\item
  The approximation satisfies the error bound $\Norm[F]{\AugTns{X} - \AugTns{C} \times_1 \AugMat{U}_1 \cdots \times_d \AugMat{U}_d} \leq \tau \Norm[F]{\AugTns{X}}$.
\item
  The core maintains a factorization of the form $\AugTns{C} = \AugTns{V} \times_d \AugMat{S}$ along the last mode, where $\AugTnsMat{V}{d}$ is row-orthogonal, and $\AugMat{S}$ is diagonal with $\AugMat{S}(i, i) > 0$ for $1 \leq i \leq \Aug{r}_d$.
\end{enumerate}
These properties exactly match those listed in Assumption~\ref{asm:streaming_setup}.
The updated factorization can then be used as is to assimilate subsequent tensor slices.

We must emphasize that, in our implementation we never form the augmented tensor $\AugTns{X}$ explicitly.
It merely serves as a convenient placeholder for analyzing the error of the approximation.

\subsection{Non-Streaming Mode Updates}
\label{sec:non-streaming-mode}

Intuitively, the incremental SVD update for matrices, described in \cref{sec:incremental-svd}, can be reframed as follows.
When adding a new row to an existing SVD factorization, we first expand the basis along the non-streaming mode (i.e., the right singular vectors) as needed so that the new data lies in their span.
Then we use the SVD of a smaller matrix to put the overall factorization in the SVD format.
In the higher-dimensional tensor setup, we have multiple non-streaming modes, so we will synchronize the basis (the factor matrices) multiple times.

Let us assume we have already synchronized the factor matrices for tensor modes $1, \ldots, k - 1$.
The current approximation of original data tensor $\Tns{X}$ consists of updated core $\Tns{C}_{k - 1} \in \RR^{\Aug{r}_1 \times \cdots \Aug{r}_{k - 1} \times r_k \times \cdots \times r_d}$, updated factor matrices $\AugMat{U}_\ell \in \RR^{n_\ell \times \Aug{r}_\ell}$ for $1 \leq \ell \leq k - 1$, and remaining factor matrices $\Mat{U}_\ell \in \RR^{n_k \times r_k}$ for $k \leq \ell \leq d$:
\begin{equation}
  \label{eq:non-stream-init-old}
  \Tns{X} \approx \Tns{C}_{k - 1} \times_1 \AugMat{U}_1 \cdots \times_{k - 1} \AugMat{U}_{k - 1} \times_k \Mat{U}_k \cdots \times_d \Mat{U}_d,
\end{equation}
and the core admits a factorization of the form
\begin{equation}
  \label{eq:non-stream-init-core}
  \Tns{C}_{k - 1} = \Tns{V}_{k - 1} \times_d \Mat{S},
\end{equation}
where $\Mat{S}$ is diagonal and the matricization $\GroupTnsMat{V}{k - 1}{d}$ has mutually orthonormal rows.
Further suppose that we have a partial factorization of the new data slice $\Tns{Y}$ obtained through successive projection onto the updated factor matrices,
\begin{align*}
  \PartCore{D}{k - 1} = \Tns{Y} \times_1 \AugMat{U}_1^\top \cdots \times_{k - 1} \AugMat{U}_{k - 1}^\top,
\end{align*}
leading to the approximation
\begin{equation}
  \label{eq:non-stream-init-new}
  \Tns{Y} \approx \PartCore{D}{k - 1} \times_1 \AugMat{U}_1 \cdots \times_{k - 1} \AugMat{U}_{k - 1}.
\end{equation}
Before the first mode have been processed, we denote $\Tns{C}_0 = \Tns{C}$, $\Tns{V}_0 = \Tns{V}$ and $\PartCore{D}{0} = \Tns{Y}$.

\begin{algorithm}[t]
  \caption{Non Streaming Mode Update}
  \label{alg:non-streaming-update}
  \begin{algorithmic}[1]
    \Function{NonStreamingUpdate}{$k, \Mat{U}_k, \Tns{C}_{k - 1}, \PartCore{D}{k - 1}, \Tns{V}_{k - 1}, \epsilon$}
      \State {$\Tns{P}_k \gets \PartCore{D}{k - 1} \times_k \Mat{U}_k^\top$}
      \State {$\Tns{E}_k \gets \PartCore{D}{k - 1} - \Tns{P}_k \times_k \Mat{U}_k$}
      \If {$\Norm[F]{\Tns{E}_k} \leq \epsilon$}
        \State {$\AugMat{U}_k \gets \Mat{U}_k$}
        \State {$\Tns{C}_k \gets \Tns{C}_{k - 1}$}
        \State {$\PartCore{D}{k} \gets \Tns{P}_k$}
        \State {$\Tns{V}_k \gets \Tns{V}_{k - 1}$}
      \Else
        \State {$r_k', \Mat{W}_k \gets \Call{TruncatedSVD}{\GroupTnsMat{E}{k}{k}, \epsilon}$}
        \State {$\Tns{L}_k \gets \Tns{E}_k \times_k \Mat{W}_k^\top$}
        \State {$\AugMat{U}_k \gets [\Mat{U}_k, \Mat{W}_k]$}
        \State {$\Tns{C}_k \gets \Call{PadWithZeros}{k, \Tns{C}_{k - 1}, r_k'}$} \Comment{See \cref{eq:pad_core_with_zeros}}
        \State {$\PartCore{D}{k} \gets \Call{Concatenate}{k, \Tns{P}_k, \Tns{L}_k}$} \Comment{See \cref{eq:update_slice_for_next_mode}}
        \State {$\Tns{V}_k \gets \Call{PadWithZeros}{k, \Tns{V}_{k - 1}, r_k'}$} \Comment{See \cref{eq:pad_singular_vectors_with_zeros}}
      \EndIf
      \State {\Return $\AugMat{U}_k, \Tns{C}_k, \PartCore{D}{k}, \Tns{V}_k$}
    \EndFunction \Comment{$\Tns{C}_k \times_k \AugMat{U}_k = \Tns{C}_{k - 1} \times_k \Mat{U}_k, \; \Norm[F]{\PartCore{D}{k - 1} - \PartCore{D}{k} \times_k \AugMat{U}_k} \leq \epsilon$}
  \end{algorithmic}
\end{algorithm}

\begin{figure}[t!]
  \centering
  \includegraphics[width=0.98\linewidth]{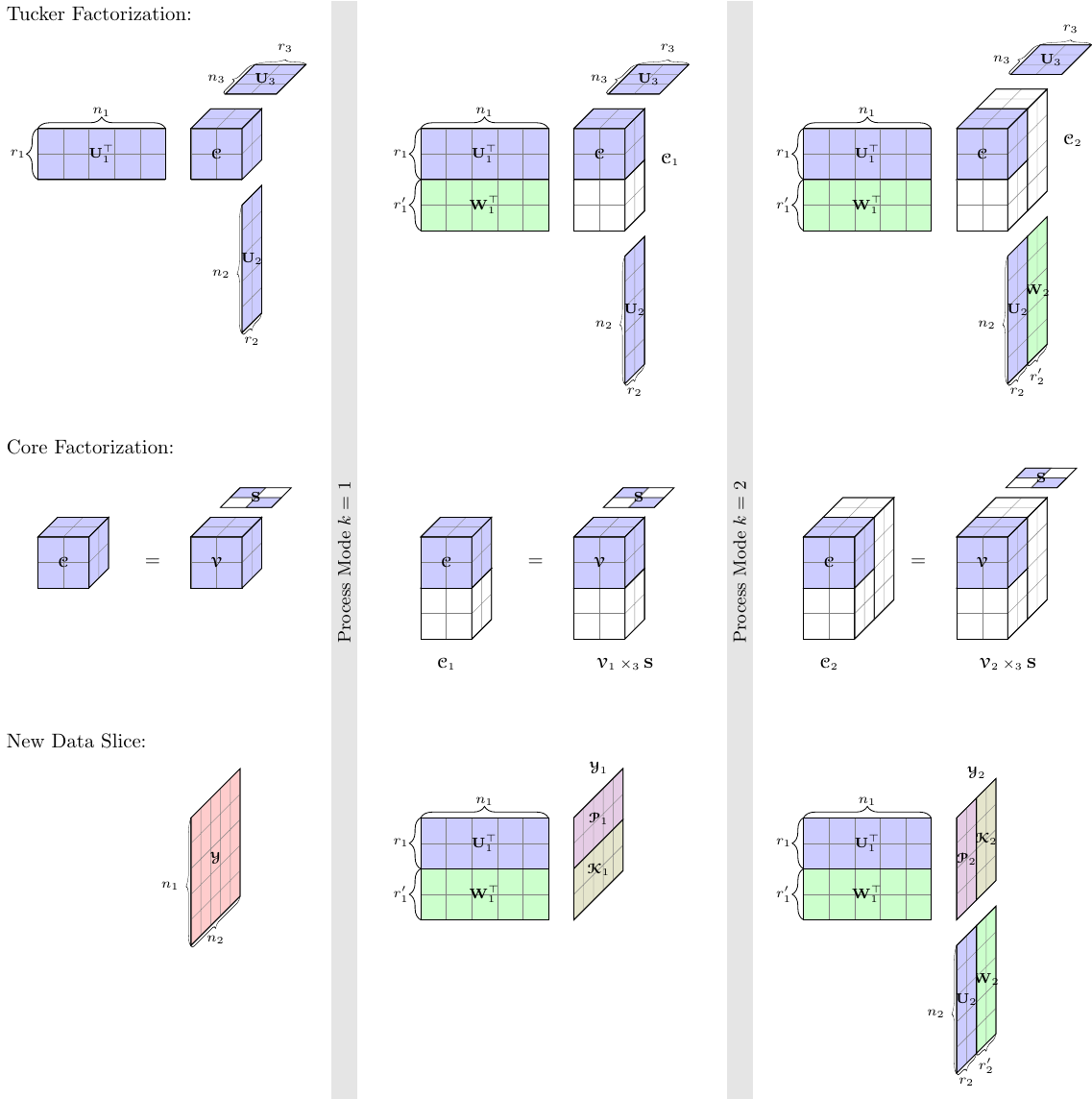}
  \caption{%
    Visualization of the streaming Tucker update steps for the non-streaming modes.
    \textit{Left:} We start from the Tucker factorization $\{\Tns{C}, \Mat{U}_1, \Mat{U}_2, \Mat{U}_3\}$ of a 3D tensor, corresponding auxiliary tensor $\Tns{V}$ needed for incremental SVD (ISVD) satisfying $\Tns{C} = \Tns{V} \times_3 \Mat{S}$ with diagonal $\Mat{S}$ with positive diagonal entries and $ \TnsMat{V}{d} \TnsMat{V}{d}^\top = \Mat{I}$, and new slice $\Tns{Y}$.
    \textit{Middle:} We process mode $k = 1$ of the new slice $\Tns{Y}$ and find orthogonal basis vectors supplementing the columns of $\Mat{U}_1$ such that $\TnsMat{Y}{1} = \Mat{U}_1 \Mat{P}_1 + \Mat{W}_1 \Mat{K}_1$. We append $\Mat{W}_1$ to update the first factor matrix $\AugMat{U}_1$, and correspondingly add zero slices to the core and auxiliary tensor to obtain $\Tns{C}_1$ and $\Tns{V}_1$; they still satisfy $\Tns{C}_1 = \Tns{V}_1 \times_3 \Mat{S}$.
    We use $\Mat{P}_1$ and $\Mat{K}_1$ to construct the projected new slice $\Tns{Y}_1$.
    \textit{Right:} We similarly process mode $k = 2$, while in principle ignoring the first factor matrix. After processing all the non-streaming modes, we are left with the projected new slice $\Tns{Y}_2 \in \RR^{(r_1 + r_1') \times (r_2 + r_2')}$, which we then insert using the ISVD procedure.%
  }
  \label{fig:streaming}
\end{figure}

For general $k$, we synchronize the mode-$k$ factor matrix as follows.
First split the new data slice $\PartCore{D}{k - 1}$ into components that are parallel and orthogonal to the column span of $\Mat{U}_k$ along the $k$-th mode:
\begin{equation*}
  \Tns{P}_k = \PartCore{D}{k - 1} \times_k \Mat{U}_k^\top, \quad \Tns{E}_k = \PartCore{D}{k - 1} - \Tns{P}_k \times_k \Mat{U}_k.
\end{equation*}
Then construct a basis $\Mat{W}_k$ consisting of $r_k' \geq 0$ vectors for the complement along the $k$-th mode, and project it along this new basis: $\Tns{L}_{k} = \Tns{E}_{k} \times_k \Mat{W}_k^\top$.
Enrich the mode-$k$ factor by adding columns of $\Mat{W}_k$: $\AugMat{U}_k = [\Mat{U}_k, \Mat{W}_k]$, and correspondingly update the core by adding $r_k'$ zero slices:
\begin{equation}
  \label{eq:pad_core_with_zeros}
  \Tns{C}_k(i_1, \ldots, i_d) =
  \begin{cases}
    \Tns{C}_{k - 1}(i_1, \ldots, i_d) & \text{if} \quad 1 \leq i_k \leq r_k, \\
    0 & \text{if} \quad r_k + 1 \leq i_k \leq \Aug{r}_k := r_k + r_k',
  \end{cases}
\end{equation}
so that $\Tns{C}_k \times_k \AugMat{U}_k = \Tns{C}_{k - 1} \times_k \Mat{U}_k$.
Next, update the new tensor slice by concatenating $\Tns{P}_k$ and $\Tns{L}_k$ along the $k$-th mode:
\begin{equation}
  \label{eq:update_slice_for_next_mode}
  \PartCore{D}{k}(i_1, \ldots, i_d) =
  \begin{cases}
    \Tns{P}_{k - 1}(i_1, \ldots, i_d), & \!\! 1 \leq i_k \leq r_k, \\
    \Tns{L}_{k - 1}(i_1, \ldots, i_{k - 1}, i_k - r_k, i_{k + 1}, \ldots, i_d), & \!\! r_k + 1 \leq i_k \leq \Aug{r}_k;
  \end{cases}
\end{equation}
this constructed tensor satisfies $\PartCore{D}{k} = \PartCore{D}{k - 1} \times \AugMat{U}_k^\top$.
Finally, to maintain the orthogonal factorization of the updated core $\Tns{C}_k$ along the last tensor mode, construct the corresponding $\Tns{V}_k$ by inserting $r_k'$ zero-slices to $\Tns{V}_{k - 1}$ along the $k$-th mode:
\begin{equation}
  \label{eq:pad_singular_vectors_with_zeros}
  \Tns{V}_k(i_1, \ldots, i_d) =
  \begin{cases}
    \Tns{V}_{k - 1}(i_1, \ldots, i_d) & \text{if} \quad 1 \leq i_k \leq r_k, \\
    0 & \text{if} \quad r_k + 1 \leq i_k \leq \Aug{r}_k := r_k + r_k'.
  \end{cases}
\end{equation}
Given \cref{eq:non-stream-init-core}, it is easy to verify that $\Tns{C}_k = \Tns{V}_k \times_d \Mat{S}$.
Furthermore, the matricization $\GroupTnsMat{V}{k}{d}$ is constructed from $\GroupTnsMat{V}{k - 1}{d}$ by inserting zero columns, which does not change the ``mutually orthonormal rows'' property.

The core steps of this update is presented in \cref{alg:non-streaming-update}, and a schematic is shown in \cref{fig:streaming}.
This construction produces updated approximations to the augmented data tensor $\AugTns{X}$ given by
\begin{align*}
  \Tns{X} &\approx \Tns{C}_k   \times_1 \AugMat{U}_1 \cdots \times_k \AugMat{U}_k \times_{k + 1} \Mat{U}_{k + 1} \cdots \times_d \Mat{U}_d, \\
  \Tns{Y} &\approx \PartCore{D}{k} \times_1 \AugMat{U}_1 \cdots \times_k \AugMat{U}_k.
\end{align*}
The approximation error over all the non-streaming mode updates is characterized as follows:

\begin{proposition}
  \label{prop:non-stream-error-bound}
  Suppose $\Tns{C}, \Mat{U}_1, \ldots, \Mat{U}_d$ constitutes the initial Tucker factorization of data tensor $\Tns{X}$, and let $\Tns{Y}$ be a new data slice along the last tensor mode.
  Then the non-streaming mode updates produce a new factorization of $\Tns{X}$ consisting of Tucker core $\Tns{C}_{d - 1}$ and factors $\AugMat{U}_1, \ldots, \AugMat{U}_{d - 1}, \Mat{U}_d$ satisfying
  \begin{equation*}
    \Norm[F]{\Tns{X} - \Tns{C}_{d - 1} \times_1 \AugMat{U}_1 \cdots \times_{d - 1} \AugMat{U}_{d - 1} \times_d \Mat{U}_d} = \Norm[F]{\Tns{X} - \Tns{C} \times_1 \Mat{U}_1 \cdots \times_d \Mat{U}_d}.
  \end{equation*}
  In addition, approximation of data slice $\Tns{Y}$ through partial cores $\PartCore{D}{0} = \Tns{Y}, \PartCore{D}{1}, \ldots, \PartCore{D}{d - 1}$ defined w.r.t.\ factor matrices $\AugMat{U}_1, \ldots, \AugMat{U}_{d - 1}$ satisfies the usual ST-HOSVD error bound
  \begin{equation*}
    \Norm[F]{\Tns{Y} - \PartCore{D}{d - 1} \times_1 \AugMat{U}_1 \cdots \times_{d - 1} \AugMat{U}_{d - 1}}^2 = \sum_{k = 1}^{d - 1} \Norm[F]{\PartCore{D}{k - 1} - \PartCore{D}{k} \times_k \AugMat{U}_k}^2.
  \end{equation*}
\end{proposition}

\begin{proof}
  See \cref{app:non-stream-error-bound}.
\end{proof}

\subsection{Updating Along Streaming Mode}
\label{sec:streaming-mode}

For the final step of the streaming update, we need to compute the principal column subspace of the unfolding matrix
\begin{equation}
  \label{eq:streaming_unfolding_update}
  \AugTnsMat{X}{d} =
  \begin{bmatrix}
    \TnsMat{X}{d} \\
    \Vec{y}^\top
  \end{bmatrix},
  \quad \Vec{y} = \text{vec}(\Tns{Y}).
\end{equation}
Using the matrix form of Tucker factorization as described in \cref{eq:tucker_matrix}, we have
\begin{equation*}
  \TnsMat{X}{d} \approx \Mat{U}_d \GroupTnsMat{C}{d - 1}{d} (\AugMat{U}_{d - 1} \otimes \cdots \otimes \AugMat{U}_1)^\top,
  \quad
  \Vec{y}^\top \approx \Vec{d}^\top (\AugMat{U}_{d - 1} \otimes \cdots \otimes \AugMat{U}_1)^\top,
\end{equation*}
where $\Vec{d} = \text{vec}(\PartCore{D}{d - 1})$.
Furthermore, the Tucker core factorization from \cref{eq:non-stream-init-core} implies
\begin{equation*}
  \GroupTnsMat{C}{d - 1}{d} = \Mat{S} \Mat{V}^\top, \quad \Mat{V} = \GroupTnsMat{V}{d - 1}{d}^\top;
\end{equation*}
the matrix $\Mat{S}$ is diagonal, and the columns of $\Mat{V}$ form a partial orthonormal basis.
Given these observations, \cref{eq:streaming_unfolding_update} simplifies to:
\begin{equation*}
  \AugTnsMat{X}{d} \approx
  \begin{bmatrix}
    \Mat{U}_d \Mat{S} \Mat{V}^\top \\
    \Vec{d}^\top
  \end{bmatrix}
  (\AugMat{U}_{d - 1} \otimes \cdots \otimes \AugMat{U}_1)^\top,
\end{equation*}
and we want to find the principal column subspace of the RHS.
This is precisely the problem solved by incremental SVD \cite{Brand02,Cheng18}, and we follow the steps outlined in \cref{sec:incremental-svd}.
Once we have obtained the principal column subspace basis $\AugMat{U}_d$ using \cref{alg:isvd}, we then update the core by demanding
\begin{equation*}
  \begin{split}
    \begin{bmatrix}
      \Mat{U}_d \GroupTnsMat{C}{d - 1}{d} \\
      \Vec{d}^\top
    \end{bmatrix}
    \approx \AugMat{U}_d \AugTnsMat{C}{d}
    \impliedby
    \AugTnsMat{C}{d}
    &= \AugMat{U}_d^\top
    \begin{bmatrix}
      \Mat{U}_d \GroupTnsMat{C}{d - 1}{d} \\
      \Vec{d}^\top
    \end{bmatrix} \\
    &= \AugMat{U}_d(1 : n_d, :)^\top \Mat{U}_d \GroupTnsMat{C}{d - 1}{d} + \AugMat{U}_d(n_d + 1, :)^\top \Vec{d}^\top.
  \end{split}
\end{equation*}
This update creates an approximation of the form $\AugTnsMat{X}{d} \approx \AugMat{U}_d \AugTnsMat{C}{d} (\AugMat{U}_{d - 1} \otimes \cdots \otimes \AugMat{U}_1)^\top$, which is precisely the Tucker factorization in matrix form.
In tensor form, the above update reads
\begin{equation*}
  \AugTns{C} = \Tns{C}_{d - 1} \times_d \AugMat{U}_d(1 : n_d, :)^\top \Mat{U}_d + \PartCore{D}{d - 1} \times_d \AugMat{U}_d(n_d + 1, :)^\top.
\end{equation*}

\subsection{The Full Algorithm}

\begin{algorithm}[t]
  \caption{Streaming ST-HOSVD Update}
  \label{alg:sthosvd_update}
  \begin{algorithmic}[1]
    \Require 
    \Statex Tucker approximation $\Tns{C}$, $\Mat{U}_1, \ldots, \Mat{U}_d$ of $\Tns{X}$ s.t.\ $\Norm[F]{\Tns{X} - \Tns{C} \times_1 \Mat{U}_1 \cdots \times_d \Mat{U}_d} \leq \tau \Norm[F]{\Tns{X}}$
    \Statex Orthogonal factorization $\Tns{C} = \Tns{V} \times_d \Mat{S}$ along the streaming mode
    \Statex New tensor slice $\Tns{Y} \in \RR^{n_1 \times \cdots \times n_{d - 1} \times 1}$
    \Ensure
    \Statex Apprx.\ $\AugTns{C}$, $\AugMat{U}_1, \ldots, \AugMat{U}_d$ of $\AugTns{X}$ defined in \cref{eq:augmented_tensor} s.t.\ $\Norm[F]{\AugTns{X} - \AugTns{C} \times_1 \AugMat{U}_1 \cdots 
    \times_d \AugMat{U}_d} \leq \tau \Norm[F]{\AugTns{X}}$
    \Statex Maintained orthogonal factorization $\AugTns{C} = \AugTns{V} \times_d \AugMat{S}$ along the streaming mode
    \Statex
    \State {Compute truncation parameter: $\epsilon \gets \tau \Norm[F]{\Tns{Y}} / \sqrt{d}$}
    \State {Initialize: $\Tns{C}_0 \gets \Tns{C}, \PartCore{D}{0} \gets \Tns{Y}, \Tns{V}_0 \gets \Tns{V}$}
    \For{$k = 1$ to $d - 1$}
      \State {$\AugMat{U}_k, \Tns{C}_k, \PartCore{D}{k}, \Tns{V}_k \gets \Call{NonStreamingUpdate}{k, \Mat{U}_k, \Tns{C}_{k - 1}, \PartCore{D}{k - 1}, \Tns{V}_{k - 1}, \epsilon}$}
    \EndFor
    \State {$\Vec{d} \gets \text{vec}(\PartCore{D}{d - 1})$}
    \State {$\Mat{V} \gets \GroupTnsMat{V}{d - 1}{d}^\top$}
    \State {$\AugMat{U}_d, \AugMat{S}, \AugMat{V} \gets \Call{IncrementalSVD}{\Mat{U}_d, \Mat{S}, \Mat{V}, \Vec{d}, \epsilon}$}
    \State {Construct $\AugTns{V}$ such that $\AugTnsMat{V}{d} = \AugMat{V}^\top$}
    \State {$\AugTns{C} \gets \Tns{C}_{d - 1} \times_d \AugMat{U}_d(1 : n_d, :)^\top \Mat{U}_d + \PartCore{D}{d - 1} \times_d \AugMat{U}_d(n_d + 1, :)^\top$}
  \end{algorithmic}
\end{algorithm}

The full streaming Tucker update scheme, along with specific approximation tolerances to ensure the appropriate error bounds, are summarized in \cref{alg:sthosvd_update}.
The justification for the error analysis is provided in \cref{app:full-algo-error-bound}.

To analyze the computational complexity of our proposed algorithm, we recap/concretize the following notations:
\begin{itemize}
  \item
    The size of the initial data tensor $\Tns{X}$ is $n_1 \times \cdots \times n_d$.
  \item
    In the initial Tucker factorization of $\Tns{X}$, the size of the Tucker core $\Tns{C}$ is $r_1 \times \cdots \times r_d$.
    Denote $r_* = r_1 \cdots r_d$.
  \item
    The size of the new data slice $\Tns{Y}$ is $n_1 \times \cdots \times n_{d - 1} \times 1$.
    For $1 \leq k \leq d - 1$, denote $n_{>k}' = n_{k + 1} \cdots n_{d - 1}$.
  \item
    After the streaming update, the updated Tucker core $\AugTns{C}$ has size $\Aug{r}_1 \times \cdots \times \Aug{r}_{d}$.
    For $1 \leq k \leq d - 1$, we define $r_k' = \Aug{r}_k - r_k$ and denote $\Aug{r}_{<k} = \Aug{r}_1 \cdots \Aug{r}_{k - 1}$.
\end{itemize}
We now quantify the computational complexity of the streaming ST-HOSVD update:

\begin{proposition}
  \label{prop:sthosvd_update_flop_count}
  The computational cost of adding a single $n_1 \times \cdots \times n_{d - 1} \times 1$ data slice to an existing $(r_1, \ldots, r_d)$-rank Tucker factorization of $n_1 \times \cdots \times n_d$ data tensor using \cref{alg:sthosvd_update} scales as
  \begin{equation*}
    \begin{split}
      O\bigg(\sum_{k = 1}^{d - 1} \Aug{r}_{<k} r_k n_k n_{>k}'
      &+ \sum_{k = 1}^{d - 1} \mathbbm{1}_{\Aug{r}_k > r_k} \left(\Aug{r}_{<k} n_k^2 n_{>k}' + n_k^3 + \Aug{r}_{<k} r_k' n_k n_{>k}' \right) \\
      &+ \Aug{r}_{<d} r_d + r_d^3 + r_d \Aug{r}_d n_d + \Aug{r}_{<d} \Aug{r}_d r_d\bigg),
    \end{split}
  \end{equation*}
  where $\Aug{r}_1, \ldots, \Aug{r}_d$ are the updated ranks, and $\mathbbm{1}_T$ is the indicator function (equals 1 if $T$ is true, equals 0 otherwise).
  The algorithm uses $O(\max\{n_{<d}, n_1^2, \ldots, n_{d - 1}^2, r_d^2, n_d \Aug{r}_d\})$ additional memory for this update step.
\end{proposition}

\begin{proof}
  See \cref{app:full-algo-flop-count}
\end{proof}

The first term of the computational cost expression corresponds to the inescapable TTM costs from the non-streaming mode updates, the second term to the cost from non-streaming mode updates that activate only when the Tucker ranks increase, and final term to the cost of updating streaming mode factorization.
This expression is obviously too complicated to draw any intuition from, so we make the simplifying assumption that Tucker ranks do not change.
In particular, if $d \geq 3$, $n_k = O(n)$ for $1 \leq k \leq d - 1$, and $r_k = O(r)$ for $1 \leq k \leq d$, then the computational complexity of inserting a single data slice simplifies to
\begin{equation*}
  O\left(\sum_{k = 1}^{d - 1} r^k n^{d - k} + r^2 n_d + r^{d + 1}\right),
\end{equation*}
and the algorithm uses additional $O(\max\{n^{d - 1}, r_d^2, n_d \Aug{r}_d\})$ memory.
Adding up the computational cost as the number of slices $n_d$ along the streaming mode grows, the amortized cost of incrementally building the Tucker factorization is $O(n_d^2)$, where we treat $n$ and $r$ as constants.
Contrast this against the batch ST-HOSVD algorithm, where the leading order cost is the eigensolve of the Gram matrix along the streaming mode: $O(n_d^3)$.

\subsection{Comparison Against Prior Works}

We provide a qualitative contrast of our approach against some of the streaming Tucker approaches mentioned earlier in \cref{sec:related_work}.
We focus on two set of works.

Sketching-based single-pass algorithms \cite{MaBe18,Sun2020} accelerate computation of the Tucker factorization by:
(a) first using structured random matrices to construct ``sketches'' that capture essential properties of the data tensor (e.g.\ span of the tensor along each mode),
(b) then constructing the Tucker factorization efficiently from these sketches.
Since sketching is a linear operation, it is highly amenable to streaming data; in fact, the data tensor does not need to arrive slice by slice, \emph{any} arbitrary order of the tensor entries will suffice.
However, every time additional data arrives, we need to repeat step (b) anew.
This may limit the applicability of the algorithms.

DTucker and DTuckerO algorithms \cite{JaKa2023} reduce the computational and memory complexity algorithms by first applying randomized SVD to 2D slices of the data tensor, then constructing/updating the Tucker factorization using these
smaller-sized intermediate factorizations.
However, these algorithms need to know \emph{a priori} the target Tucker ranks of the factorization.

Both of these approaches rely on randomized linear algebra algorithms; consequently, any theoretical result provides probabilistic inequalities on the error (e.g. bounds on expected squared error).
In contrast, our proposed algorithm is entirely deterministic, and we always end up with a fully realized Tucker factorization at the end of each streaming update step.
Our update strategy ensures that the relative error between observed data and the Tucker compressed representation remains smaller than a user-specified relative tolerance set at the beginning of the streaming compression.
Additionally, the algorithm allows the Tucker ranks to grow as required, eliminating the need to \emph{a priori} set the Tucker ranks; however, any rank growth is accompanied by extra computations that in turn adds to the runtime.


\section{Results}
\label{sec:results}

We implemented our streaming ST-HOSVD algorithm within the TuckerMPI C++ software package\footnote{\url{https://gitlab.com/tensors/TuckerMPI}}.
Note that the implementation we are testing is serial and largely leverages the kernels (e.g. TTM, Gram computation etc.) that are already part of TuckerMPI.
We compare our streaming algorithm against the standard batch ST-HOSVD on synthetic and simulation data tensors.
All tests were run on a single node of a compute cluster equipped with two 16-core, 2.10 GHz Intel Xeon E5-2683 v4 CPUs and 256 GB of RAM.

The experiments follow the same general pattern.
Specifically, given a data tensor $\Tns{X}$ and a target tolerance $\tau$, we use first the batch ST-HOSVD algorithm to construct the reference Tucker-compressed representation and record the obtained ranks, runtime and memory utilization.
Next, we apply our streaming ST-HOSVD algorithm to compress the data tensor (starting the computation with a fraction of all slices), and compare the obtained ranks, runtime,
and memory utilization against those from the batch version.
To assess the accuracy of the streaming algorithm, we reconstruct the tensor from the both batch and streaming Tucker compressed representation, and compute the relative reconstruction error $\epsilon = \Norm[F]{\Tns{X} - \widehat{\Tns{X}}} / \Norm[F]{\Tns{X}}$, where $\widehat{\Tns{X}}$ is the reconstructed tensor from the compressed representations.

\subsection{Synthetic dataset}

Consider the $d$-dimensional tensor $\Tns{F} \in \RR^{n_1 \times \cdots \times n_d}$, constructed from sampling the function
\begin{equation*}
  f(x_1, \ldots, x_d) = \sum_{j_1 = -f_1}^{f_1} \cdots \sum_{j_d = -f_d}^{f_d} c_{j_1, \ldots, j_d} \sin{(j_1 x_1 + \cdots + j_d x_d) \pi}
\end{equation*}
on a uniform grid in the domain $[-1, 1]^d$. Using trigonometric ``sum of angles'' identities, we can easily establish that the exact rank of the $k$-th unfolding $\Mat{F}_{(k)}$ is $r_k = 2 f_k + 1$, and the $d$-dimensional Tucker core of this tensor has dimensions $r_1 \times \cdots \times r_d$.

We construct two $100 \times 100 \times 2500$ tensors with ranks $(11, 11, 11)$ and $(11, 11, 15)$ and concatenate them; this combined tensor $\overline{\Tns{X}} \in \RR^{100 \times 100 \times 5000}$ thus undergoes a ``rank growth'' halfway through the streaming mode.
We define our data tensor to be $\Tns{X} = \overline{\Tns{X}} + \Tns{N}$ by adding noise; the entries of the additive noise tensor are drawn from the standard Gaussian distribution and scaled such that $\Norm[F]{\Tns{N}} = \eta \Norm[F]{\Tns{X}}$ for a noise-to-signal ratio $\eta = 10^{-3}$.
We compress this ``sine wave'' data tensor using batch ST-HOSVD and streaming ST-HOSVD with error tolerance $\tau = 10^{-2}$.

\begin{figure}[t]
  \centering
  \begin{subfigure}{0.48\linewidth}
    \centering
    \includegraphics[width=\linewidth]{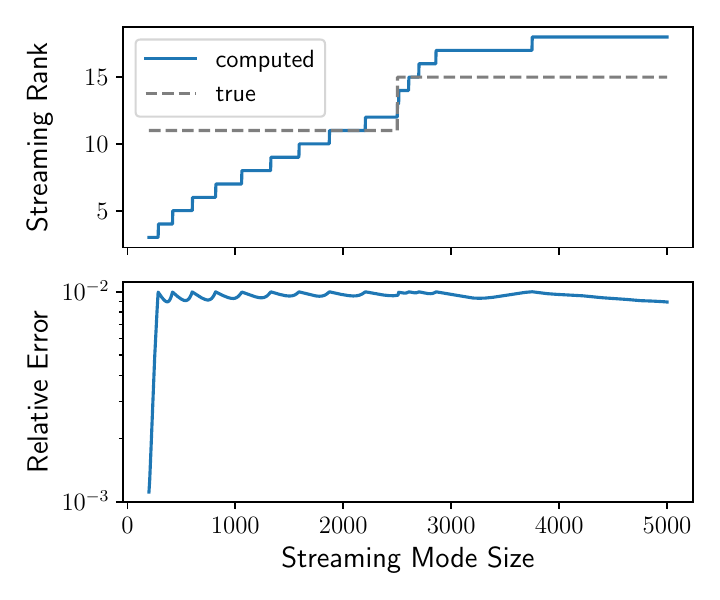}
  \end{subfigure}
  ~
  \begin{subfigure}{0.48\linewidth}
    \centering
    \includegraphics[width=\linewidth]{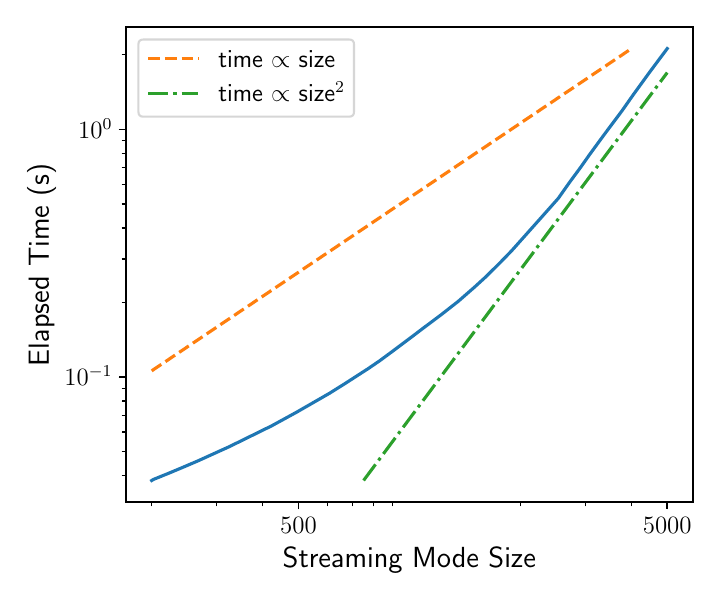}
  \end{subfigure}
  \caption{%
    Streaming compression of the $100 \times 100 \times 5000$ sine-wave dataset which has approximate Tucker rank $(11, 11, 11)$ up to the 2500-th slice, but whose overall ranks are $(11, 11, 15)$.
    \emph{Left:} The growth of the Tucker rank $r_3$ along the streaming mode, and evolution of the relative reconstruction error, as we insert slices one by one using our streaming algorithm after the initial ST-HOSVD with 200 slices.
    \emph{Right:} Quadratic scaling of the streaming Tucker algorithm runtime w.r.t.\ streaming mode size.
  }
  \label{fig:noisy_wave}
\end{figure}

In \cref{fig:noisy_wave}, we plot the variation of the Tucker ranks and the reconstruction error as we construct the streaming approximation starting from an initial ST-HOSVD with 200 slices.
We note that our streaming algorithm detects the rank bump at the 2500\textsuperscript{th} slice and adjusts the streaming rank accordingly while maintaining the prescribed error bound.
We also observe that the runtime of the streaming algorithm grows quadratically with the streaming mode size, as opposed to the theoretical cubic growth of the batch algorithm runtime.

\begin{table}[t]
  \centering
  \caption{%
    Comparison of batch and streaming ST-HOSVD methods with the $100 \times 100 \times 5000$ sine-wave dataset with approximate Tucker ranks $(11, 11, 15)$.
    The streaming algorithm is initialized from an ST-HOSVD with 200 snapshots.
    We report the Tucker ranks, maximum ST-HOSVD memeory usage, and computation time.%
  }
  \label{tbl:noisy_wave}
  \medskip
  \begin{tabular}{c c c c}
    \toprule
    Algorithm &  Rank          & Memory (MB) & Time (s) \\
    \midrule
    Batch     & $(11, 11, 15)$ &       605.7 &    194.6 \\
    Stream    & $(11, 11, 18)$ &        18.3 &      2.1 \\
    \bottomrule
  \end{tabular}
\end{table}

We compare the final performance of the streaming algorithm against the standard batch ST-HOSVD in \cref{tbl:noisy_wave}; the results indicate both the streaming and batch versions recover comparable Tucker ranks from the datasets.
However, we see an advantage in terms of the memory usage and computational time.
During the incremental updates, our algorithm uses far less memory compared to the batch version.
Moreover, the streaming algorithm is much faster compared to the batch algorithm for this dataset, where the streaming rank does not increase drastically.
This reduces the overall computational time when the streaming mode size is larger than the non-streaming mode sizes.

\subsection{HCCI combustion dataset}

Following a previous study \cite{KollaKC2020}, we consider a simulation of auto-ignition processes of a turbulent ethanol-air mixture, under conditions corresponding to a homogeneous charge compression ignition (HCCI) combustion mode of an internal combustion engine, to assess the streaming ST-HOSVD algorithm applied to scientific data.
The simulation was performed on a rectangular Cartesian grid of a 2D spatial domain with $672 \times 672$ grid points, which constitute the first two modes of the tensor.
The simulation state comprises 33 variables at each grid point---the third mode---and 268 time snapshots (the last of a total of 626 timesteps) are considered which constitutes the streaming mode of a fourth-order tensor.
Since this is not a synthetic data set with an exact known low-rank structure, the results of the streaming and batch ST-HOSVD algorithms are likely to produce different results due to our lack of knowledge about the exact truncation criteria; however they should still be comparable for a fixed relative error tolerance $\tau$.

\begin{figure}[t]
  \centering

  \begin{subfigure}{0.31\textwidth}
    \centering
    \includegraphics[width=\linewidth]{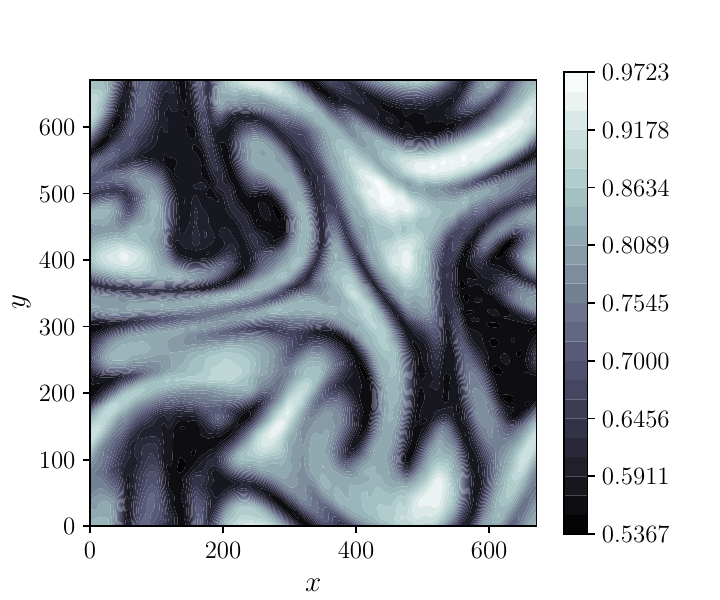}
    \caption{Compressed, $\tau = 10^{-1}$}
  \end{subfigure}
  ~
  \begin{subfigure}{0.31\textwidth}
    \centering
    \includegraphics[width=\linewidth]{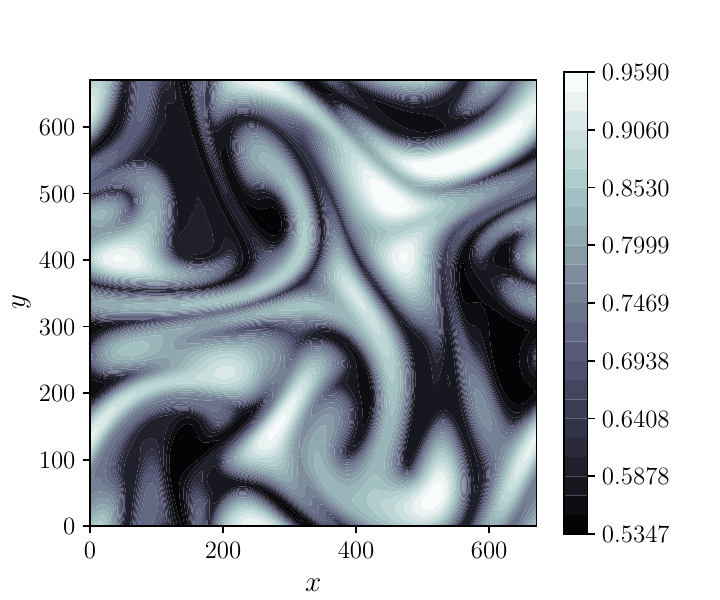}
    \caption{Compressed, $\tau = 10^{-2}$}
  \end{subfigure}
  ~
  \begin{subfigure}{0.31\textwidth}
    \centering
    \includegraphics[width=\linewidth]{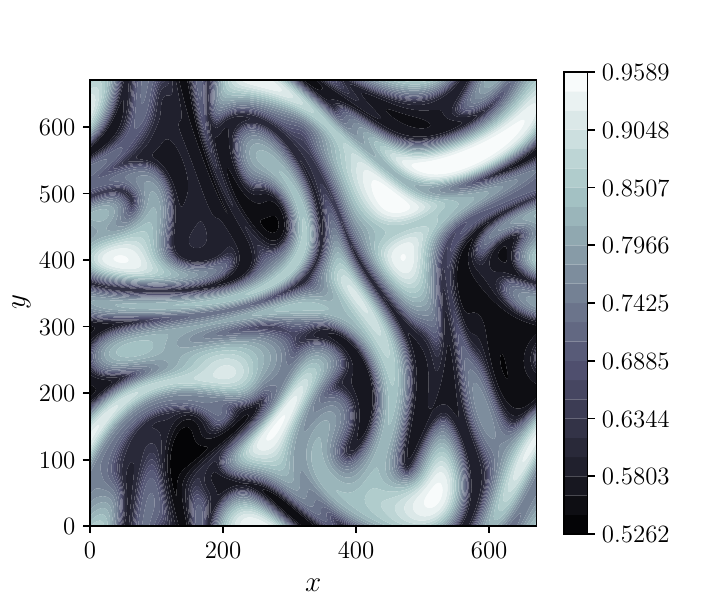}
    \caption{Original}
  \end{subfigure}

  \caption{%
    Visual comparison between the original data tensor and reconstructed streaming Tucker tensor at compression tolerances $\tau \in \{10^{-1}, 10^{-1}\}$.
    We consider the temperature variable across the spatial modes at time $t = 2 \times 10^{-3}$~s, and reconstruct it from the streaming Tucker compressed tensor initialized with 100 snapshots.%
  }
  \label{fig:hcci}
\end{figure}

In \cref{fig:hcci}, we compare the reconstructions of the temperature variable at simulation time $t = 2.0 \times 10^{-3}$~s from the streaming Tucker compression (initialized with 100 timesteps of interest) against the original dataset at two tolerance levels: $\tau \in \{10^{-1}, 10^{-2}\}$.
We observe that, as the tolerance is tightened, the contours of the element-wise reconstructions become visually nearly identical to the contours of the original data.

\begin{table}[t]
  \centering
  \caption{%
    Performance of the batch and streaming ST-HOSVD methods when compressing the $672 \times 672 \times 33 \times 268$ HCCI combustion dataset at two ST-HOSVD tolerances.
    Streaming compressions are initialized with 100 and 150 snapshots.
    We report the Tucker ranks, maximum ST-HOSVD memeory usage (in GBs), computation time (in seconds), and estimated relative reconstruction errors.%
  }
  \label{tbl:hcci}
  \medskip
  \begin{tabular}{c c c c c c}
    \toprule
    Tolerance & Algorithm    & Rank                 & Memory & Time   & Relative Error       \\
    \midrule
    $10^{-1}$ & Batch        & $( 35,  27,  8,  7)$ & 33.68  & 119.48 & $9.6 \times 10^{-2}$ \\
              & Stream (100) & $( 46,  38, 10, 12)$ & 12.78  &  86.23 & $9.4 \times 10^{-2}$ \\
              & Stream (150) & $( 43,  34,  9, 10)$ & 19.09  &  95.96 & $9.4 \times 10^{-1}$ \\
    \midrule                                                 
    $10^{-2}$ & Batch        & $(122, 110, 21, 29)$ & 38.70  & 152.86 & $9.6 \times 10^{-3}$ \\
              & Stream (100) & $(140, 128, 24, 44)$ & 14.88  & 197.68 & $9.9 \times 10^{-3}$ \\
              & Stream (150) & $(134, 121, 23, 37)$ & 22.09  & 168.83 & $9.9 \times 10^{-3}$ \\
    \bottomrule
  \end{tabular}
\end{table}

In \cref{tbl:hcci}, we compare the performance of the batch and streaming Tucker factorizations.
We compress the data tensor with two error tolerances, and initialize the streaming versions with data from the first 100 and 150 timesteps of interest.
We observe that the target error tolerance is satisfied in all the experiments.
Additionally, as we increase the number of initial time steps, the Tucker ranks from the streaming algorithm approach those computed from the batch algorithm.
This improvement, however, comes at the cost of increasing memory footprint of the streaming algorithm (which is still smaller compared to the batch version).

These experiments with the HCCI dataset demonstrate that our streaming ST-HOSVD implementation can process moderately large datasets (tens of gigabytes) on single node compute clusters.
Compressing larger datasets (terabytes and higher) will however require a MPI-based fully distributed implementation \cite{DeKoMePhRi2024} to overcome the limitations on per-node memory availability.


\section{Conclusions and Future Work}
\label{sec:conclusions}

In this paper, we introduced an algorithm for constructing Tucker factorizations of streaming data and analyzed its theoretical properties (computational and memory complexity, approximation error bounds).
Starting from an initial ST-HOSVD factorization, the proposed algorithm updates the compressed representation as new tensor slices are incorporated one at a time while ensuring the relative reconstruction error stays within a specified fixed tolerance.
In our experiments, we observed that the streaming version may on occasion overestimate the Tucker ranks when compared to the batch algorithm, but the former generally has a smaller memory footprint.
Additionally, this rank-overestimation can be controlled by supplying a larger number of snapshots to the initial ST-HOSVD.
Finally, when the size of the streaming mode is large compared to the product of Tucker ranks along the non-streaming modes, the streaming version usually outperforms the batch algorithm in terms of computational time.
In summary, our algorithm is particularly suited to approximately low-rank ($r_k \ll n_k$) tensors with large streaming mode size ($n_d \gg r_1 \cdots r_{d - 1}$).

The current version of the algorithm processes a single slice at a time and uses a basic version of incremental SVD.
In future, we want to enable processing multiple data slices at once, and utilize the advanced incremental SVD algorithms developed in \cite{Brand02,Cheng18}.
This will likely reduce dependency on the size of initial batch factorization, stabilize the rank growth, and lead to better approximations of data tensors.

\bibliographystyle{plain}
\bibliography{references}

\appendix
\newpage


\section{Proof of \Cref{prop:isvd-error-bound}: Incremental SVD Error Bound}
\label{app:isvd-error-bound}

Consider the $r$-rank approximate SVD of $\Mat{A} \in \RR^{m \times n}$,
\begin{equation*}
  \Mat{A} = \Mat{U} \Mat{S} \Mat{V}^\top + \Mat{E},
\end{equation*}
where $\Mat{U}$ and $\Mat{V}$ are orthogonal matrices satisfying $\Mat{U}^\top \Mat{U} = \Mat{V}^\top \Mat{V} = \Mat{I}_r$.
Given a new row $\Vec{b} \in \RR^{n}$, the incremental SVD update is conditioned on whether we ignore the orthogonal complement $\Vec{e} = \Vec{b} - \Vec{V} \Vec{p}$, where $\Vec{p} = \Mat{V}^\top \Vec{b}$.
If ignored, the approximation error is given by
\begingroup
\allowdisplaybreaks
\begin{align*}
  \AugMat{A} - \AugMat{U} \AugMat{S} \AugMat{V}^\top
  &=
  \begin{bmatrix}
    \Mat{A} \\
    \Vec{b}^\top
  \end{bmatrix}
  -
  \begin{bmatrix}
    \Mat{U} &   \\
          & 1
  \end{bmatrix}
  \Mat{U}'' \AugMat{S} \Mat{V}^{\prime\prime,\top}
  \Mat{V}^\top \\
  &=
  \begin{bmatrix}
    \Mat{A} \\
    \Vec{b}^\top
  \end{bmatrix}
  -
  \begin{bmatrix}
    \Mat{U} &   \\
          & 1
  \end{bmatrix}
  \begin{bmatrix}
    \Mat{S} \\
    \Vec{p}^\top
  \end{bmatrix}
  \Mat{V}^\top \\
  &=
  \begin{bmatrix}
    \Mat{A} \\
    \Vec{b}^\top
  \end{bmatrix}
  -
  \begin{bmatrix}
    \Mat{U} \Mat{S} \Mat{V}^\top \\
    \Vec{p}^\top \Mat{V}^\top
  \end{bmatrix} \\
  &=
  \begin{bmatrix}
    \Mat{A} - \Mat{U} \Mat{S} \Mat{V}^\top \\
    (\Vec{b} - \Mat{V} \Vec{p})^\top
  \end{bmatrix} \\
  &=
  \begin{bmatrix}
    \Mat{A} - \Mat{U} \Mat{S} \Mat{V}^\top \\
    \Vec{e}^\top
  \end{bmatrix},
\end{align*}
\endgroup
so that $\Norm[F]{\AugMat{A} - \AugMat{U} \AugMat{S} \AugMat{V}^\top}^2 = \Norm[F]{\Mat{A} - \Mat{U} \Mat{S} \Mat{V}^\top}^2 + \Norm{\Vec{e}}^2$.
On the other hand, if we construct a new basis $\Vec{v} = q^{-1} \Vec{e}$ the orthogonal complement where $q = \Norm{\Vec{e}}$, then $\Vec{b} = \Mat{V} \Vec{p} + q \Vec{v}$ and the approximation error is given by
\begingroup
\allowdisplaybreaks
\begin{align*}
  \AugMat{A} - \AugMat{U} \AugMat{S} \AugMat{V}^\top
  &=
  \begin{bmatrix}
    \Mat{A} \\
    \Vec{b}^\top
  \end{bmatrix}
  -
  \begin{bmatrix}
    \Mat{U} &   \\
          & 1
  \end{bmatrix}
  \Mat{U}'' \AugMat{S} \Mat{V}^{\prime\prime,\top}
  \Mat{V}^{\prime,\top} \\
  &=
  \begin{bmatrix}
    \Mat{A} \\
    \Vec{b}^\top
  \end{bmatrix}
  -
  \begin{bmatrix}
    \Mat{U} &   \\
          & 1
  \end{bmatrix}
  \begin{bmatrix}
    \Mat{S}      &   \\
    \Vec{p}^\top & q
  \end{bmatrix}
  \begin{bmatrix}
    \Mat{V}^\top \\
    \Vec{v}^\top
  \end{bmatrix} \\
  &=
  \begin{bmatrix}
    \Mat{A} \\
    \Vec{b}^\top
  \end{bmatrix}
  -
  \begin{bmatrix}
    \Mat{U} \Mat{S} \Mat{V}^\top \\
    \Vec{p}^\top \Mat{V}^\top + q \Vec{v}^\top
  \end{bmatrix} \\
  &=
  \begin{bmatrix}
    \Mat{A} - \Mat{U} \Mat{S} \Mat{V}^\top \\
    (\Vec{b} - \Mat{V} \Vec{p} - q \Vec{v})^\top
  \end{bmatrix} \\
  &=
  \begin{bmatrix}
    \Mat{A} - \Mat{U} \Mat{S} \Mat{V}^\top \\
    \null
  \end{bmatrix},
\end{align*}
\endgroup
so that $\Norm[F]{\AugMat{A} - \AugMat{U} \AugMat{S} \AugMat{V}^\top}^2 = \Norm[F]{\Mat{A} - \Mat{U} \Mat{S} \Mat{V}^\top}^2$.
This concludes our proof.

\section{Proof of \Cref{prop:isvd-flop-count}: Incremental SVD Complexity}
\label{app:isvd-flop-count}

We itemize the number of floating point operations in incremental SVD update of an $r$-rank SVD of a $m \times n$ data matrix:
\begin{itemize}
\item
  $O(n r)$ operations and $O(n)$ extra storage in computing the projection $\Vec{p}$ of the new row $\Vec{b}$ onto the row span of $\Mat{V}$.
\item
  $O(n r)$ operations and $O(n)$ extra storage to compute the orthogonal complement $\Vec{e}$ of the new row w.r.t.\ the row span of $\Mat{V}$.
\item
  $O(n)$ operations and $O(1)$ extra storage to determine the norm $q$ of $\Vec{e}$, and if needed, another $O(n)$ operations and $O(n)$ extra storage to determine the additional basis vector $\Vec{v}$.
\item
  $O(n r^2)$ operations $O(n r)$ storage to form the matrix $\Mat{V}'$ of size $n \times r$ or $n \times (r + 1)$ and optionally compute its QR factorization for orthogonalization.
\item
  $O((r + 1)^3) = O(r^3)$ operations and $O((r + 1)^2) = O(r^2)$ extra storage to form the matrix $\Mat{S}'$ of size either $(r + 1) \times (r + 1)$ or $(r + 1) \times r$, followed by its SVD to compute the factors $\Mat{U}''$, $\AugMat{S}$ and $\Mat{V}''$.
  Assume that this decomposition reveals rank $\Aug{r}$, i.e. $\AugMat{S} \in \RR^{\Aug{r} \times \Aug{r}}$; then $\Aug{r} \leq r + 1$.
\item
  $O((m + 1) (r + 1) \Aug{r}) = O(m r^2)$ and $O(n (r + 1) \Aug{r}) = O(n r^2)$ operations, and $O((m + 1) \Aug{r}) = O(m r)$ and $O(n \Aug{r}) = O(nr)$ extra storage, to compute the $\AugMat{U}$ and $\AugMat{V}$ via matrix multiplication.
\end{itemize}
Summing all these contributions, we end up with an algorithm with $O(r^3 + m r^2 + n r^2)$ operations and $O(r^2 + m r + n r)$ extra storage in the worst case.

To compute the amortized cost when the rank $r$ is fixed, we sum up contributions from individual incremental SVD steps.
Assuming we start from a data matrix with $m_\text{init}$ rows, the total number of floating point operations is
\begin{equation*}
  \begin{split}
    O\left(\sum_{\ell = m_\text{init}}^{m - 1} r^3 + \ell r^2 + n r^2\right)
    &= O\left(\sum_{\ell = 0}^{m - 1} r^3 + \ell r^2 + n r^2\right) \\
    &= O\left(m r^3 + \frac{(m - 1) m}{2} r^2 + m n r^2\right) \\
    &= O(m r^3 + m^2 r^2 + m n r^2).
  \end{split}
\end{equation*}
The amortized storage complexity is simply the maximum over all incremental steps: $O(r^2 + m r + n r)$.

\section{Proof of \Cref{prop:non-stream-error-bound}: Non-Streaming Update Error Bound}
\label{app:non-stream-error-bound}

Given the two sets of factor matrices $\{\AugMat{U}_1, \ldots, \AugMat{U}_{k - 1}\}$ and $\{\Mat{U}_{k + 1}, \ldots, \Mat{U}_d\}$, define
\begin{align*}
  \AugMat{U}_{<k} &= \Mat{I}_{r_d} \otimes \cdots \otimes \Mat{I}_{r_{k + 1}} \otimes \AugMat{U}_{k - 1} \otimes \cdots \otimes \AugMat{U}_1, \\
  \Mat{U}_{>k} &= \Mat{U}_d \otimes \cdots \otimes \Mat{U}_{k + 1} \otimes \Mat{I}_{n_{k - 1}} \otimes \cdots \otimes \Mat{I}_{n_1}.
\end{align*}
By the mixed-product property of the Kronecker product and using the matricized form of TTMs, we can rewrite \cref{eq:non-stream-init-old,eq:non-stream-init-new} as
\begin{alignat*}{3}
  & \TnsMat{X}{k} &&\approx \Mat{U}_k \GroupTnsMat{C}{k - 1}{k} (\Mat{U}_d \otimes \cdots \otimes \Mat{U}_{k + 1} \otimes \AugMat{U}_{k - 1} \otimes \cdots \otimes \AugMat{U}_1)^\top &&= \Mat{U}_k \GroupTnsMat{C}{k - 1}{k} \Mat{U}_{>k}^\top \AugMat{U}_{<k}^\top, \\
  & \TnsMat{Y}{k} &&\approx \PartCoreMat{D}{k - 1}{k} (\Mat{I}_{r_d} \otimes \cdots \otimes \Mat{I}_{r_{k + 1}} \otimes \AugMat{U}_{k - 1} \otimes \cdots \otimes \AugMat{U}_1)^\top  &&= \PartCoreMat{D}{k - 1}{k} \AugMat{U}_{<k}^\top.
\end{alignat*}
Then the mode-$k$ matricization of the augmented tensor in \cref{eq:augmented_tensor} is given by
\begin{equation*}
  \begin{split}
    \AugTnsMat{X}{k}
    =
    \begin{bmatrix}
      \TnsMat{X}{k} & \TnsMat{Y}{k}
    \end{bmatrix}
    &\approx
    \begin{bmatrix}
      \Mat{U}_k \GroupTnsMat{C}{k - 1}{k} \Mat{U}_{>k}^\top \AugMat{U}_{<k}^\top & \PartCoreMat{D}{k - 1}{k} \AugMat{U}_{<k}^\top
    \end{bmatrix} \\
    &=
    \begin{bmatrix}
      \Mat{U}_k \GroupTnsMat{C}{k - 1}{k} \Mat{U}_{>k}^\top & \PartCoreMat{D}{k - 1}{k}
    \end{bmatrix}
    \AugMat{U}_{<k}^\top,
  \end{split}
\end{equation*}
Given the definitions of $\Tns{P}_k$, $\Tns{E}_k$ and $\Tns{L}_k$, we can relate their mode-$k$ matricizations as
\begin{alignat*}{4}
  & \Tns{P}_k &&= \PartCore{D}{k - 1} \times_k \Mat{U}_k^\top        \quad &&\iff \quad \GroupTnsMat{P}{k}{k} &&= \Mat{U}_k^\top \PartCoreMat{D}{k - 1}{k}, \\
  & \Tns{E}_k &&= \PartCore{D}{k - 1} - \Tns{P}_k \times_k \Mat{U}_k \quad &&\iff \quad \GroupTnsMat{E}{k}{k} &&= \PartCoreMat{D}{k - 1}{k} - \Mat{U}_k \GroupTnsMat{P}{k}{k}, \\
  & \Tns{L}_k &&= \Tns{E}_k \times_k \Mat{W}_k^\top                  \quad &&\iff \quad \GroupTnsMat{L}{k}{k} &&= \Mat{W}_k^\top \GroupTnsMat{E}{k}{k}.
\end{alignat*}
Then it follows that
\begin{equation*}
  \GroupTnsMat{E}{k}{k} \approx \Mat{W}_k \GroupTnsMat{L}{k}{k} \implies \PartCoreMat{D}{k - 1}{k} \approx \Mat{U}_k \GroupTnsMat{P}{k}{k} + \Mat{W}_k \GroupTnsMat{L}{k}{k},
\end{equation*}
and consequently,
\begin{align*}
  \begin{split}
    \AugTnsMat{X}{k}
    &\approx
    \begin{bmatrix}
      \Mat{U}_k \GroupTnsMat{C}{k - 1}{k} \Mat{U}_{>k}^\top & \Mat{U}_k \GroupTnsMat{P}{k}{k} + \Mat{W}_k \GroupTnsMat{L}{k}{k}
    \end{bmatrix}
    \AugMat{U}_{<k}^\top \\
    &= 
    \begin{bmatrix}
      \Mat{U}_k & \Mat{W}_k
    \end{bmatrix}
    \begin{bmatrix}
      \GroupTnsMat{C}{k - 1}{k} \Mat{U}_{>k}^\top & \GroupTnsMat{P}{k}{k} \\
                                             & \GroupTnsMat{L}{k}{k}
    \end{bmatrix}
    \AugMat{U}_{<k}^\top \\
    &= \AugMat{U}_k
    \begin{bmatrix}
      \GroupTnsMat{C}{k}{k} \Mat{U}_{>k}^\top & \PartCoreMat{D}{k}{k}
    \end{bmatrix}
    \AugMat{U}_{<k}^\top,
  \end{split}
\end{align*}
where the last step follows from the definitions of $\AugMat{U}_k$, $\Tns{C}_k$ and $\PartCore{D}{k}$.
Thus, we have updated the approximations of $\Tns{X}$ and $\Tns{Y}$ as
\begin{alignat*}{4}
  & \TnsMat{X}{k} &&\approx \AugMat{U}_k \GroupTnsMat{C}{k}{k} \Mat{U}_{>k}^\top \AugMat{U}_{<k}^\top \quad &&\iff \quad \Tns{X} &&\approx \Tns{C}_k \times_1 \AugMat{U}_1 \cdots \times_k \AugMat{U}_k \times_{k + 1} \Mat{U}_{k + 1} \cdots \times_d \Mat{U}_d, \\
  & \TnsMat{Y}{k} &&\approx \AugMat{U}_k \PartCoreMat{D}{k}{k} \AugMat{U}_{<k}^\top                   \quad &&\iff \quad \Tns{Y} &&\approx \PartCore{D}{k} \times_1 \AugMat{U}_1 \cdots \times_k \AugMat{U}_k.
\end{alignat*}
Now, since $\Mat{U}_k$ and $\Mat{W}_k$ have orthonormal column spans, we have
\begin{align*}
  \Mat{W}_k^\top \PartCoreMat{D}{k - 1}{k}
  &= \Mat{W}_k^\top (\Mat{U}_k \GroupTnsMat{P}{k}{k} + \GroupTnsMat{E}{k}{k}) \\
  &= \Mat{W}_k^\top \Mat{U}_k \GroupTnsMat{P}{k}{k} + \Mat{W}_k^\top \GroupTnsMat{E}{k}{k} \\
  &= \Mat{0} \cdot \GroupTnsMat{P}{k}{k} + \GroupTnsMat{L}{k}{k} \\
  &= \GroupTnsMat{L}{k}{k}.
\end{align*}
It follows that,
\begin{equation*}
  \AugMat{U}_k \GroupTnsMat{C}{k}{k}
  =
  \begin{bmatrix}
    \Mat{U}_k & \Mat{W}_k
  \end{bmatrix}
  \begin{bmatrix}
    \GroupTnsMat{C}{k - 1}{k} \\
    \null
  \end{bmatrix}
  = \Mat{U}_k \GroupTnsMat{C}{k - 1}{k},
\end{equation*}
and
\begin{equation*}
  \AugMat{U}_k^\top \PartCoreMat{D}{k - 1}{k}
  =
  \begin{bmatrix}
    \Mat{U}_k^\top \\
    \Mat{W}_k^\top
  \end{bmatrix}
  \PartCoreMat{D}{k - 1}{k}
  =
  \begin{bmatrix}
    \Mat{U}_k^\top \PartCoreMat{D}{k - 1}{k} \\
    \Mat{W}_k^\top \PartCoreMat{D}{k - 1}{k} \\
  \end{bmatrix}
  =
  \begin{bmatrix}
    \GroupTnsMat{P}{k}{k} \\
    \GroupTnsMat{L}{k}{k}
  \end{bmatrix}
  = \PartCoreMat{D}{k}{k};
\end{equation*}
in tensor format, these read $\Tns{C}_k \times_k \AugMat{U}_k = \Tns{C}_{k - 1} \times_k \Mat{U}_k$ and $\PartCore{D}{k - 1} \times_k \AugMat{U}_k^\top = \PartCore{D}{k}$, respectively.

Starting from $\Tns{C}_0 = \Tns{C}$ and $\PartCore{D}{0} = \Tns{Y}$ and combing these relations over modes $1$ through $k$, we obtain:
\begin{equation*}
  \begin{split}
    \Tns{C} \times_1 \Mat{U}_1 \cdots \times_k \Mat{U}_k
    &= \Tns{C}_0 \times_1 \Mat{U}_1 \cdots \times_k \Mat{U}_k \\
    &= \Tns{C}_1 \times_1 \AugMat{U}_1 \times_2 \Mat{U}_2 \cdots \times_k \Mat{U}_k \\
    &= \Tns{C}_2 \times_1 \AugMat{U}_1 \times_2 \AugMat{U}_2 \times_3 \Mat{U}_3 \cdots \times_k \Mat{U}_k \\
    &= \cdots \\
    &= \Tns{C}_k \times_1 \AugMat{U}_1 \cdots \times_k \AugMat{U}_k,
  \end{split}
\end{equation*}
and
\begin{equation*}
  \begin{split}
    \PartCore{D}{k}
    &= \PartCore{D}{k - 1} \times_k \AugMat{U}_k \\
    &= \PartCore{D}{k - 2} \times_{k - 1} \AugMat{U}_{k - 1} \times_k \AugMat{U}_k \\
    &= \cdots \\
    &= \PartCore{D}{0} \times_1 \AugMat{U}_1 \cdots \times_k \AugMat{U}_k \\
    &= \Tns{Y} \times_1 \AugMat{U}_1 \cdots \times_k \AugMat{U}_k.
  \end{split}
\end{equation*}
Given this successive approximation of the new data slice $\Tns{Y}$ through partial cores in the style of ST-HOSVD, we establish the error bound
\begin{equation*}
  \Norm[F]{\Tns{Y} - \PartCore{D}{k} \times_1 \AugMat{U}_1 \cdots \times_k \AugMat{U}_k}^2 = \sum_{\ell = 1}^{k} \Norm[F]{\PartCore{D}{\ell - 1} - \PartCore{D}{\ell} \times_\ell \AugMat{U}_\ell}^2.
\end{equation*}
Taking $k = d - 1$ completes the proof.

\section{Full Algorithm Error Analysis}
\label{app:full-algo-error-bound}

We start by splitting the approximation error of a Tucker approximation $\Tns{X} \approx \Tns{C} \times_1 \Mat{U}_1 \cdots \times_d \Mat{U}_d$ along the non-streaming and streaming modes:
\begin{align*}
  \tau \Norm[F]{\Tns{X}}^2
  &\geq \Norm[F]{\Tns{X} - \Tns{C} \times_1 \Mat{U}_1 \cdots \times_d \Mat{U}_d}^2 \\
  &= \Norm[F]{\TnsMat{X}{d} - \Mat{U}_d \TnsMat{C}{d} \Mat{U}_{<d}^\top}^2 \\
  &= \Norm[F]{\TnsMat{X}{d} - \TnsMat{X}{d} \Mat{U}_{<d} \Mat{U}_{<d}^\top + (\TnsMat{X}{d} \Mat{U}_{<d} - \Mat{U}_d \TnsMat{C}{d}) \Mat{U}_{<d}^\top}^2 \\
  &= \Norm[F]{\TnsMat{X}{d} - \TnsMat{X}{d} \Mat{U}_{<d} \Mat{U}_{<d}^\top}^2 + \Norm[F]{\TnsMat{X}{d} \Mat{U}_{<d} - \Mat{U}_d \TnsMat{C}{d}}^2.
\end{align*}
A similar split for the approximation error after the streaming Tucker update yields
\begin{align*}
  \Norm[F]{\AugTns{X} - \AugTns{C} \times_1 \AugMat{U}_1 \cdots \times_d \AugMat{U}_d}^2
  &= \Norm[F]{\AugTnsMat{X}{d} - \AugMat{U}_d \AugTnsMat{C}{d} \AugMat{U}_{<d}^\top}^2 \\
  &= \Norm[F]{\AugTnsMat{X}{d} - \AugTnsMat{X}{d} \AugMat{U}_{<d} \AugMat{U}_{<d}^\top}^2 + \Norm[F]{\AugTnsMat{X}{d} \AugMat{U}_{<d} - \AugMat{U}_d \AugTnsMat{C}{d}}^2.
\end{align*}
Now, after the non-streaming mode updates, we obtain the approximation to the new slice,
\begin{equation*}
  \Tns{D} = \Tns{Y} \times_1 \AugMat{U}_1^\top \cdots \times_{d - 1} \AugMat{U}_{d - 1}^\top
  \implies
  \Vec{d}^\top = \Vec{y}^\top \AugMat{U}_{<d},
\end{equation*}
where $\Vec{y}$ and $\Vec{d}$ are vectorizations of $\Tns{Y}$ and $\Tns{D}$, respectively.
It follows
\begin{equation*}
  \AugTnsMat{X}{d} \AugMat{U}_{<d}
  =
  \begin{bmatrix}
    \TnsMat{X}{d} \\
    \Vec{y}^\top
  \end{bmatrix}
  \AugMat{U}_{<d}
  =
  \begin{bmatrix}
    \TnsMat{X}{d} \AugMat{U}_{<d} \\
    \Vec{y}^\top \AugMat{U}_{<d}
  \end{bmatrix}
  =
  \begin{bmatrix}
    \TnsMat{X}{d} \AugMat{U}_{<d} \\
    \Vec{d}^\top
  \end{bmatrix},
\end{equation*}
and consequently
\begin{align*}
  &\phantom{\null=\null} \Norm[F]{\AugTns{X} - \AugTns{C} \times_1 \AugMat{U}_1 \cdots \times_d \AugMat{U}_d}^2 \\
  &= \Norm[F]{\AugTnsMat{X}{d} - \AugTnsMat{X}{d} \AugMat{U}_{<d} \AugMat{U}_{<d}^\top}^2 + \Norm[F]{\AugTnsMat{X}{d} \AugMat{U}_{<d} - \AugMat{U}_d \AugTnsMat{C}{d}}^2 \\
  &= \left\lVert
    \begin{bmatrix}
      \TnsMat{X}{d} \\
      \Vec{y}^\top
    \end{bmatrix}
    -
    \begin{bmatrix}
      \TnsMat{X}{d} \AugMat{U}_{<d} \AugMat{U}_{<d}^\top \\
      \Vec{d}^\top \AugMat{U}_{<d}^\top
    \end{bmatrix}
  \right\rVert_F^2
  +
  \left\lVert
    \begin{bmatrix}
      \TnsMat{X}{d} \AugMat{U}_{<d} \\
      \Vec{d}^\top
    \end{bmatrix}
    -
    \AugMat{U}_d \AugTnsMat{C}{d}
  \right\rVert_F^2 \\
  &= \Norm[F]{\TnsMat{X}{d} - \TnsMat{X}{d} \AugMat{U}_{<d} \AugMat{U}_{<d}^\top}^2 + \Norm[F]{\Vec{y}^\top - \Vec{d}^\top \AugMat{U}_{<d}^\top}^2
  +
  \left\lVert
    \begin{bmatrix}
      \TnsMat{X}{d} \AugMat{U}_{<d} \\
      \Vec{d}^\top
    \end{bmatrix}
    -
    \AugMat{U}_d \AugTnsMat{C}{d}
  \right\rVert_F^2.
\end{align*}
To simplify the last term, we note that $\operatorname{cols}(\Mat{U}_k) \subseteq \operatorname{cols}(\AugMat{U}_k)$ for all $1 \leq k \leq d - 1$ by construction of the non-streaming mode updates.
It follows that $\operatorname{cols}(\Mat{U}_{<d}) \subseteq \operatorname{cols}(\AugMat{U}_{<d})$; let $\Mat{W}_{<d}$ be the matrix that consists of the columns that are in $\AugMat{U}_{<d}$, but not in $\Mat{U}_{<d}$.
Then
\begin{equation*}
  \AugMat{U}_{<d} =
  \begin{bmatrix}
    \Mat{U}_{<d} & \Mat{W}_{<d}
  \end{bmatrix}
  \Mat{\Pi},
\end{equation*}
where $\Mat{\Pi}$ is a permutation matrix inducing the appropriate column permutation required for this equality.
Denoting $\Mat{E} := \TnsMat{X}{d} \Mat{U}_{<d} - \Mat{U}_d \TnsMat{C}{d}$, we can write
\begin{align*}
  \TnsMat{X}{d} \AugMat{U}_{<d}
  &= \TnsMat{X}{d} \begin{bmatrix} \Mat{U}_{<d} & \Mat{0} \end{bmatrix} \Mat{\Pi} + \TnsMat{X}{d} \begin{bmatrix} \Mat{0} & \Mat{W}_{<d} \end{bmatrix} \Mat{\Pi} \\
  &= \begin{bmatrix} \TnsMat{X}{d} \Mat{U}_{<d} & \Mat{0} \end{bmatrix} \Mat{\Pi} + \begin{bmatrix} \Mat{0} & \TnsMat{X}{d} \Mat{W}_{<d} \end{bmatrix} \Mat{\Pi} \\
  &= \begin{bmatrix} \Mat{U}_d \TnsMat{C}{d} + \Mat{E} & \Mat{0} \end{bmatrix} \Mat{\Pi} + \begin{bmatrix} \Mat{0} & \TnsMat{X}{d} \Mat{W}_{<d} \end{bmatrix} \Mat{\Pi} \\
  &= \begin{bmatrix} \Mat{U}_d \TnsMat{C}{d} & \Mat{0} \end{bmatrix} \Mat{\Pi} + \begin{bmatrix} \Mat{E} & \Mat{0} \end{bmatrix} \Mat{\Pi} + \begin{bmatrix} \Mat{0} & \TnsMat{X}{d} \Mat{W}_{<d} \end{bmatrix} \Mat{\Pi} \\
  &= \Mat{U}_d \begin{bmatrix} \TnsMat{C}{d} & \Mat{0} \end{bmatrix} \Mat{\Pi} + \begin{bmatrix} \Mat{E} & \TnsMat{X}{d} \Mat{W}_{<d} \end{bmatrix} \Mat{\Pi} \\
  &= \Mat{U}_d \GroupTnsMat{C}{d - 1}{d} + \begin{bmatrix} \Mat{E} & \TnsMat{X}{d} \Mat{W}_{<d} \end{bmatrix} \Mat{\Pi} \\
  &= \Mat{U}_d \Mat{S} \GroupTnsMat{V}{d - 1}{d} + \begin{bmatrix} \Mat{E} & \TnsMat{X}{d} \Mat{W}_{<d} \end{bmatrix} \Mat{\Pi}
\end{align*}
Applying incremental SVD update with tolerance $\epsilon = \frac{\tau}{\sqrt{d}} \Norm[F]{\Tns{Y}}$ bounds the approximation error of the last term:
\begin{align*}
  \left\lVert
    \begin{bmatrix}
      \TnsMat{X}{d} \AugMat{U}_{<d} \\
      \Vec{d}^\top
    \end{bmatrix}
    -
    \AugMat{U}_d \AugTnsMat{C}{d}
  \right\rVert_F^2
  &\leq \left\lVert\begin{bmatrix} \Mat{E} & \TnsMat{X}{d} \Mat{W}_{<d} \end{bmatrix} \Mat{\Pi}\right\rVert_F^2 + \frac{\tau^2}{d} \Norm[F]{\Tns{Y}}^2 \\
  &= \left\lVert\begin{bmatrix} \Mat{E} & \TnsMat{X}{d} \Mat{W}_{<d} \end{bmatrix}\right\rVert_F^2 + \frac{\tau^2}{d} \Norm[F]{\Tns{Y}}^2 \\
  &= \Norm[F]{\Mat{E}}^2 + \Norm[F]{\TnsMat{X}{d} \Mat{W}_{<d}}^2 + \frac{\tau^2}{d} \Norm[F]{\Tns{Y}}^2 \\
  &= \Norm[F]{\TnsMat{X}{d} \Mat{U}_{<d} - \Mat{U}_d \TnsMat{C}{d}}^2 + \Norm[F]{\TnsMat{X}{d} \Mat{W}_{<d}}^2 + \frac{\tau^2}{d} \Norm[F]{\Tns{Y}}^2,
\end{align*}
Now, let $\Mat{Q}$ be an orthogonal matrix such that columns of $\AugMat{U}_{<d}$ and $\Mat{Q}$ form a complete orthonormal basis.
Recalling that columns of $\AugMat{U}_{<d}$ are simply union of columns of $\Mat{U}_{<d}$ and $\Mat{W}_{<d}$, it follows that columns of orthogonal matrices $\Mat{U}_{<d}$, $\Mat{W}_{<d}$ and $\Mat{Q}$ form a full orthonormal basis as well.
Thus,
\begin{align*}
  &\phantom{\null=\null} \Norm[F]{\TnsMat{X}{d} - \TnsMat{X}{d} \AugMat{U}_{<d} \AugMat{U}_{<d}^\top}^2 + \Norm[F]{\TnsMat{X}{d} \Mat{W}_{<d}}^2 \\
  &= \Norm[F]{\TnsMat{X}{d} - \TnsMat{X}{d} \AugMat{U}_{<d} \AugMat{U}_{<d}^\top}^2 + \Norm[F]{\TnsMat{X}{d} \Mat{W}_{<d} \Mat{W}_{<d}^\top}^2 \\
  &= \Norm[F]{\TnsMat{X}{d} \Mat{Q} \Mat{Q}^\top}^2 + \Norm[F]{\TnsMat{X}{d} \Mat{W}_{<d} \Mat{W}_{<d}^\top}^2 \\
  &= \left\lVert\TnsMat{X}{d} \begin{bmatrix} \Mat{W}_{<d} & \Mat{Q} \end{bmatrix} \begin{bmatrix} \Mat{W}_{<d} & \Mat{Q} \end{bmatrix}^\top \right\rVert_F^2 \\
  &= \Norm[F]{\TnsMat{X}{d} - \TnsMat{X}{d} \Mat{U}_{<d} \Mat{U}_{<d}^\top}^2,
\end{align*}
and we obtain the final error bound,
\begin{align*}
  &\phantom{\null\leq\null} \Norm[F]{\AugTns{X} - \AugTns{C} \times_1 \AugMat{U}_1 \cdots \times_d \AugMat{U}_d}^2 \\
  &\leq \Norm[F]{\TnsMat{X}{d} - \TnsMat{X}{d} \AugMat{U}_{<d} \AugMat{U}_{<d}^\top}^2 + \Norm[F]{\Vec{y}^\top - \Vec{d}^\top \AugMat{U}_{<d}^\top}^2 + \Norm[F]{\TnsMat{X}{d} \Mat{U}_{<d} - \Mat{U}_d \TnsMat{C}{d}}^2 + \null \\
  &\qquad\qquad \Norm[F]{\TnsMat{X}{d} \Mat{W}_{<d}}^2 + \frac{\tau^2}{d} \Norm[F]{\Tns{Y}}^2 \\
  &= \Norm[F]{\TnsMat{X}{d} - \TnsMat{X}{d} \AugMat{U}_{<d} \AugMat{U}_{<d}^\top}^2 + \Norm[F]{\TnsMat{X}{d} \Mat{W}_{<d}}^2 + \Norm[F]{\TnsMat{X}{d} \Mat{U}_{<d} - \Mat{U}_d \TnsMat{C}{d}}^2 + \null \\
  &\qquad\qquad \Norm[F]{\Vec{y}^\top - \Vec{d}^\top \AugMat{U}_{<d}^\top}^2 + \frac{\tau^2}{d} \Norm[F]{\Tns{Y}}^2 \\
  &= \Norm[F]{\TnsMat{X}{d} - \TnsMat{X}{d} \Mat{U}_{<d} \Mat{U}_{<d}^\top}^2 + \Norm[F]{\TnsMat{X}{d} \Mat{U}_{<d} - \Mat{U}_d \TnsMat{C}{d}}^2 + \Norm[F]{\Vec{y}^\top - \Vec{d}^\top \AugMat{U}_{<d}^\top}^2 + \frac{\tau^2}{d} \Norm[F]{\Tns{Y}}^2 \\
  &\leq \tau \Norm[F]{\Tns{X}}^2 + \frac{d - 1}{d} \tau^2 \Norm[F]{\Tns{Y}}^2 + \frac{\tau^2}{d} \Norm[F]{\Tns{Y}}^2 \\
  &= \tau^2 (\Norm[F]{\Tns{X}}^2 + \Norm[F]{\Tns{Y}}^2) \\
  &= \tau^2 \Norm[F]{\AugTns{X}}^2.
\end{align*}

\section{Proof of \cref{prop:sthosvd_update_flop_count}: Full Algorithm Complexity}
\label{app:full-algo-flop-count}

The most expensive parts of \cref{alg:sthosvd_update} are the tensor-times-matrix products, Gram computation and the eigendecomposition; we will focus on these operations.

In updating the $k$-th non-streaming mode using \cref{alg:non-streaming-update}, we carry out at minimum two TTMs.
Computational complexity of both of these are identical: $O(\Aug{r}_{<k} r_k n_k n_{>k}')$.
In the worst case, where the existing basis $\Mat{U}_k$ is not sufficient to capture the new data slice (leading to $\Aug{r}_k > r_k$), we need to carry out a Gram computation, an eigendecomposition, and a TTM.
The cost of these three operations are $O(\Aug{r}_{<k} n_k^2 n_{>k})$, $O(n_k^3)$, and $O(\Aug{r}_{<k} r_k' n_k n_{>k}')$, respectively.

For the streaming mode update using \cref{alg:isvd}, the major costs stem from the matrix-vector products needed for projection, truncated SVD, and TTMs.
These operations cost $O(\Aug{r}_{<d} r_d)$, $O(r_d^3)$, and $O(r_d \Aug{r}_d n_d + \Aug{r}_{<d} \Aug{r}_d r_d)$, respectively.
Combining the leading order costs of the non-streaming and streaming mode updates, we obtain the flop count in \cref{prop:sthosvd_update_flop_count}.

To analyze the memory complexity, we proceed as follows.
Similar to the batch ST-HOSVD algorithm, the size of the data tensor decays by a factor of $\frac{\Aug{r}_k}{n_k}$ as we process each tensor mode $k$.
Thus, the highest memory utilization from the TTMs corresponds to processing mode $k = 1$; each subsequent mode can simply re-use the auxiliary memory allocated when processing the first mode.
If the Tucker ranks increase during the non-streaming mode updates, then we also need to store $O(n_k^2)$ sized Gram matrices and their eigendecompositions.
When processing the streaming mode, we need to carry our SVD of a $O(r_d)$-sized square matrix, and temporarily store the $O(n_d \Aug{r}_d)$-sized updated factor matrix.
This leads to the memory complexity in \cref{prop:sthosvd_update_flop_count}.

\end{document}